\documentclass[preprint]{imsart}
\RequirePackage[OT1]{fontenc}
\RequirePackage{amsthm,amsmath}
\RequirePackage[numbers]{natbib}
\RequirePackage[colorlinks,citecolor=blue,urlcolor=blue]{hyperref}
 \RequirePackage{graphicx}
\usepackage{amssymb}
\usepackage{amsmath}
\usepackage{bm}
\usepackage{enumerate} 
\arxiv{arXiv:0000.0000}

\startlocaldefs
\numberwithin{equation}{section}
\theoremstyle{plain}

\endlocaldefs

\newtheorem{theorem}{Theorem}
\newtheorem{Corollary}{Corollary}
\newtheorem{Lemma}{Lemma}
\newtheorem{proposition}{Proposition}

\newtheorem{rem}{Remark}
\newtheorem{exmp}{Example}
\newtheorem{ass}{Assumption}

\def\eq{\equiv} \def\s{\sigma}

\def\beq{\begin{equation}} \def\enq{\end{equation}}

\def\e{\epsilon} \def\f{\frac}

\def\eq{\equiv}

\def\p{\partial}

\begin{document}

\begin{frontmatter}
\title{ Escape probabilities of   compound renewal processes  with drift\thanksref{}}
\runtitle{Escape probabilities}
 
\begin{aug}
\author{\fnms{Javier} \snm{Villarroel}\thanksref{t1}\ead[label=e1]{javier@usal.es}},
\author{\fnms{Juan A.} \snm{Vega}\thanksref{}\ead[label=e2]{jantovc@usal.es}}
\and
\author{\fnms{Miquel} \snm{Montero}\thanksref{t1,t2}
\ead[label=e3]{miquel.montero@ub.edu}
\ead[label=u1,url]{}}
 
\thankstext{t1}{ Spanish Agencia Estatal de Investigaci\'on FIS2016-78904 and European Fondo Europeo de Desarrollo Regional  (AEI/FEDER, UE) C3-2-P}
\thankstext{t2}{ Ag\`encia de Gesti\'o d'Ajuts Universitaris i de Recerca (AGAUR),  2017SGR1064.}
\runauthor{J. Villarroel et al.}

\affiliation{Universidad de Salamanca \thanksmark{m1}and Universitat de Barcelona \thanksmark{m2}}

\address{ Instit. Univ. de F\'{\i}sica  y Matem\'aticas  and   Dept.  Estad\'istica, \\ Plaza Merced s/n, E-37008 Salamanca, Spain \\ 
\printead{e1}
\phantom{E-mail:\ }\printead*{e2}}

\address{ Departament de F\'{\i}sica de Mat\`eria Condensada,  \\ 
Univ. de Barcelona (UB) , Mart\'{\i} i Franqu\`es 1, \\E-08028 Barcelona, Spain
\printead{e3} }
\end{aug}

\begin{abstract} We consider the  problem  of determining  escape probabilities  from an interval of   a   general compound renewal  process    with
 drift.  This problem is reduced to the solution of a certain integral equation. In  an   actuarial situation where only negative jumps arise we give a general solution for  escape   and  survival   probabilities  under Erlang$(n)$ and hypo-exponential arrivals. These   ideas are generalized to  the class of arrival distributions having rational Laplace transforms. In a general situation with two-sided jumps we also identify important families of solvable cases.      A parallelism with the ``scale function'' of 
 diffusion processes is drawn. 
  
\end{abstract}

\begin{keyword}[class=MSC]
\kwd[Primary ]{60G55}
\kwd{ 60K05}
\kwd[; secondary ]{60J75}
\end{keyword}

\begin{keyword}
\kwd{ Escape and ruin probabilities, renewal reward process}
\kwd{Integral equations}
\end{keyword}

\end{frontmatter}

\section{Introduction}   

The    problem  of determining   escape probabilities  from an interval of     general diffusions    is a classical issue  solved in terms of the ``scale and speed functions''  (see  \cite{Feller1,Bhat,kt81} for an overview).  Unfortunately,  no such   well established theory   exists  for  compound renewal processes    with
 drift. Concretely  we consider  a
random process $(X_t)_{t\ge   0} $ on a probability space $(\Omega,\mathcal G,\Bbb P)$   whose dynamics combines  uniform motion  with  speed $c\in\Bbb R$ and  sudden jumps $ J_n\in\Bbb R$  at time  epochs $ \mathbf{t}_n, n\in\Bbb N    $  triggered by a renewal process  $(N_t)_{t\ge t_0}$,  where     $  N_t= \#\{\mathbf{t_n}\in \Bbb T:  0<\mathbf{t_n}\le t\} $ counts  the number of  ``events'' $ \mathbf{t}_n, n=1,\dots\infty   $  ``observed'' in
 the time   window $ (0,t]$ and       $\Bbb T=\{\mathbf{t}_0, \mathbf{t}_1,\dots \mathbf{t}_n,\dots\}\subset\Bbb R^+$.    We define $ \mathbf{t}_0=0$ and  $x=X_{t_0}$.  Thus    
\begin{equation}
X_t=x+ ct +Y_{t}, \quad Y_t= \sum_{n=0}^{N_t}  J_n, \label{process}
\end{equation}

   When    $c=0$  the resulting ``renewal reward process''  $ X_t=\sum_{n=0}^{N_t}  J_n $ has  a prominent role      in reliability and system maintenance.
It  also describes earthquake  shocks \cite{Sornette} or stock markets  where sudden price changes  are allowed, \cite{Merton}.

   The  prototype model   of risk theory       to describe the cash flow   $(X_t)_{t\ge t_0}  $  at an insurance company results when a drift  $c>0$  is incorporated to account for the  constant   premium's     rate.    By contrast claims    arrive according to a  renewal reward process $ Y_t= \sum_{n=0}^{N_t}  J_n    $  with arrivals $\mathbf{t}_n, \mathbf{t}_n<\mathbf{t}_{n+1}$ and  sizes (or ''severities'') $J_n <0$.       This   classical {\it    risk reserve process} was first introduced by  Cramer-Lundberg under  Poissonian arrivals $N_t\sim\mathcal P(\lambda t)$ \cite{Cramer} and generalized  to general renewals   by Sparre-Andersen, cf.  \cite{Sparre}. It is usually complemented with  the  ``net profit condition'' 
$    c\Bbb E \tau_1+\Bbb E J_1  >0$ (NPC) $-$see \cite{Mikosch, Rolski} for general background.
 Even such  simplified  situation is far from trivial  and during the  last two decades  substantial research has been devoted to this topic: 
    Ruin probabilities with Poisson  arrivals  have  been studied in \cite{Cramer}.  
   Under    Erlang $\Gamma(2,\lambda)$  arrivals  they  can be  represented  as a compound geometric random variable, cf. \cite{dh0,dh01}. See also \cite{Garcia,Kluppelberg2}. The distribution of the time to ruin  under    Erlang times  is considered in\cite{dh012,Berger2,lg04,dh08}.   Ruin  probabilities under more general settings like L\'evy  and stable  processes appear in   the interesting papers 
 \cite{Kluppelberg3,Mikosch2}. See also\cite{Huzak,Kluppelberg,Wang,Dang,Lefevre}.
 However far less  is known about two-barrier exit  probabilities   even under the    assumption $cJ_1<0$. 
      
In this work       such    one-sided jump restriction  $cJ_j<0$  is not  required. To the best of our knowledge little is known about  such general  models   even though they       occur   naturally in other physical contexts: Energy dissipation in defective nonlinear
optical fibers is described by  \eqref{process} where  $c>0$ and $J_j>0$ account  for Energy  losses   due to damping  and    inhomogeneities  respectively(see \cite{MV3}).  
      Further motivation is given by the    description of  temporally
aggregated rainfall  in meteorology and hydrology contexts, see \cite{Gole}.           Here $(X_t)$ measures  rainfall accumulated at a dam    with  $J_n\ge 0 $ representing rainfall intensity   from the $j-$th shower while     a term $c t\in\Bbb R$ accounts for     the overall constant  water    inflow   rate  due to the  opposite  effects of evaporation,     water consumption  and  melting of ice and inflow of  water (hence both cases $c>0,<0$ might appear).     In a  different  context,  \eqref{process} also models 
   the  dynamics of snow depth on mountain hillsides, see \cite{Perona2}. Here  both   positive and negative  jumps may occur due to   snowfalls  and, respectively,  avalanches.   The drift term $c<0$ accounts for   snow melting    during no-snow days.  Finally, the    space dynamics   of  bacteria and several other living organisms is described by a renewal  process with a linear drift term, see \cite{Othmer}.  Exit times are  naturally related to the question  of whether  certain levels will be attained.

This paper    is structured as follows.   
   Let $a<x_0<b$ be two fixed levels and call   $\tau^a=\inf\{t>0: X_t \le   a\} $,  $\tau^b=\inf\{t>0: X_t \ge b\}>0 $.    We study {\it  two-barrier  escape  probabilities $  \Bbb P^x\Big(  \tau^b<\tau^a \Big)$, the probability  that starting from $x\in (a,b)$ the process \eqref{process}  exits  $(a,b)$ via the upper barrier}, when  both   positive and negative  jumps occur.     By means of  renewal arguments  we 
 show   that  {\it the basic EP solves  a certain   linear Fredholm integral equation (IE)}  with non-constant coefficients, cf.  eq. \eqref{Fred}.      (We use $   \Bbb P^x(A)\equiv    \Bbb P(A|X_0=x), A\in\mathcal F_\infty$). We note that {\it for pure jump Markov processes}  escape probabilities  (EP) have been  considered  by   extension of  Feller ideas and  Dynkin's formula. Some ideas in this regard appear in \cite{Oks,Duan,Liao}.  However the difficulty of the  resulting   Dirichlet problem   has prevented much progress for the solution (nevertheless, in a   remarkable paper Bertoin (\cite{Bertoin}) considers  exit   probabilities for  one sided  (i.e. without   positive jumps) stable L\'evy processes). The formalism  of Feller-Markov  semigroups is not generally  applicable here since  {\it  \eqref{process} is not Markov} (nevertheless such theory  is briefly used).   Such  lack of Markovianess implies  that  relevant  probabilities    depend on   the accessible information.    We also study how   {\it accumulated information affects more general EP }  of the form $   \Bbb P\Big(\tau^b<\tau^0  \mid
  \mathcal F_r\Big)  $ where $\mathcal F_r\equiv  \s(X_t, t\le  r) $ is the  information $\s-$ field and $r$  is an arbitrary epoch of time.  

   Once established that all EP are   codified in terms of the solution of a Fredholm IE we  devote our interest     to obtaining solutions  for the previous IE. Unfortunately, in a general situation    {\it  a closed form solution  is not possible.} Thus we attempt to classify the variety of cases that may arise (see table 1) and       clarify the role of different jump contributions. For ample classes of data we derive simplified equations   and  give the corresponding solution   (sections 4-6).  
   Due to its importance, a great deal of interest is devoted     to     the case  where  support  $(cJ_1)  \subset(-\infty,0)-$  the risk model.      Under    Poisson arrivals we give (section 4) a general solution for   the EP. We find the factorization, see \eqref{escapeexp} 
 $ \Bbb P^x\Big(\tau^b<\tau^a \Big)=\pi(x-a)/\pi(b-a)     $ 
for a certain $\pi:\Bbb R^+\to \Bbb R$ that depends on the jump distribution. $\pi$ can be identified with   the survival probability $S(x)$  when it  exists. 
  
     The analysis is then extended to  hypo-exponential arrivals  (sums of $n$ independent  exponential   variables  with  different  rates) and hence in particular to   Erlang arrivals $\tau_1\sim\Gamma(n,\lambda)$. We show that     $   \Bbb P^x\Big(\tau^b<\tau^0 \Big)=\frac{ \Delta(x,b)}{\Delta(b,b)}$ with  $  \Delta(x,b)= \det{\bm \Theta}(x,b)  $ for a certain   matrix $({\bm \Theta})_{n\times n}$.

In section 5  we develop a formalism to  deal with an ample  class of arrival distributions.
 We prove that  if $\tau_1\sim F$  has rational Laplace transform and    support  $J_1\subset(-\infty,0) $ one can derive   an equivalent   integro-differential equation,   amenable to Laplace transform.  We discuss  how to incorporate  appropriate   boundary conditions at $x=b$  that  pin down the EP.  The previous representation still holds  with a far more complicated matrix ${\bm \Theta}_{ }$.   
These ideas generalize to a double barrier situation  previous studies regarding the  survival probability $-$which could be  recovered  letting $b\to\infty$.

  Section 6 considers the problem of solving  the corresponding IE under a   general situation where jumps can take both signs. We identify  important cases  where such task can be accomplished:

When        support  $J_1\subset (-\infty,-b)\cup [0,\infty)$ the EP can be determined in closed form, {\it regardless the distribution of arrival times and jump sizes}.  

The case    support  $J_1\subset (-\infty,-0)\cup [b-x,\infty)$ is also  remarkable:  the  solution is a natural generalization of the risk model  of sections 4,5.    Finally under  Poissonian   arrivals and   jumps     with rational  characteristic function       {\it        the EP satisfies a simple ordinary differential equation}.  In  Table 1  we   summarize  our results for the EP. 
Note how scale functions only appear in  the very particular case  when  $(X)$ is  L\'evy   with negative  jumps.
   Ideas for mean escape times off $(a,b)$  appear  in \cite{MV3,MV}. 

The appendices are devoted to establish several facts  that codify densities in terms of    differential equations. We assume familiarity with  Schwartz tempered distribution theory   and  Banach's fixed point theorems.   

Given $X_0=x$ then  $t_b=(b-x)/c$  is the time  remaining  to reach the boundary $b$ when   no jumps happen. Besides  $ A^c$ is the  complementary event of $A\in\Omega$.   Given the  cdf  $F$, $\bar F=1-F$ denotes its tail; if     $g:\Bbb R\to \Bbb R$  is Borel measurable and and $\mu_F$  the Lebesgue-Stieltjes measure    associated to $F$ we set  $ \int_{\Bbb R} g(x)dF(x):=\int_{\Bbb R} g(x)d\mu_F(x)$.  
The Laplace transform (LT) of $g:\Bbb R^+\to \Bbb R $    is denoted as $  \mathcal L(g)(s)\eq \hat g(s):= \int_{0}^\infty  g(x)e^{-sx}dx$.

  \begin{table*}
\caption{\rm   The EP  $N_b(x)$ in terms of the decomposition \eqref{H} and  the characteristics of $F$ and $H$.  Recall that   the classes $ \mathcal M  $ and $\mathcal H$ are     defined in \eqref{Lapl} and  \eqref{Four}. Symbol   $  \star $ denotes a non-null component. We also list  the  IE appropriate  to particular cases and the  solution.   
}
 \begin{tabular}{@{}lllllllllll}
  \hline $q_1$ & $ q_2$   & $ p_1$& $p_2$ &   Jumps  &  $F$& $H$ & IE  & $N_b(x)$ 
&$N_b(x)$ \\   \hline 
 
 \\    $  \star $  & $  \star $  & 0 & $0 $     &   One-sided   &     $  \mathcal E(\lambda)$ & $-$ &\eqref{Fredac2} & \eqref{escape2}  &$ \pi(x)/ \pi(b)$  
 \\    $  \star $  & $  \star $  & 0 & $0 $     &   One-sided   &  Hypoexp.& $-$   &\eqref{Qeq} & \eqref{escape4} &$ \Delta(x,b)/\Delta(b,b)$ 
 \\   
   $  \star $  & $  \star $  & 0 & $0 $     & One-sided &   $  \in\mathcal M$   &     $-$&   \eqref{diffprob} &   \eqref{escape6} &$ \Delta(x,b)/\Delta(b,b)$ 
 \\   
     $  \star $  & $  \star $  & 0 & $\star $     & Two-sided &    $  \in\mathcal M$  &     $-$&   \eqref{Qeq2} &   \eqref{escape5}&$ \Delta(x,b)/\Delta(b,b)$
 \\ 
  $  0 $  & $  \star $  & 0 & $\star $     &   Two-sided &  $-$ & $-$  &       \eqref{fav2}& \eqref{double} &  $\pi(b-x)$ \\   
  $  0 $  & $  \star $  & $\star$ & $\star $     &    Two-sided& $-$  &  $-$ &  \eqref{Q2}  &  \eqref{Solfav2}&$\pi(b-x)$
\\
  $\star$   & $  \star $  & $  \star $  & $\star$    &Two-sided & $   \mathcal E(\lambda)$ &      &  \eqref{Poissongen}  &  $-$   
 \\
  $\star$   & $  \star $  & $  \star $  & $\star$    &Two-sided & $   \mathcal E(\lambda)$ &    $\in \mathcal H$   &   & \eqref{solvable2}  &$ \Delta(x,b)/\Delta(b,b)$ 
\\ 
  $  \star  $  & $  \star $  & $\star$ & $\star $     &    Two-sided& $-$  &  $-$ &\eqref{Fred}    & $-$

    \end{tabular}
\end{table*}

 \section{Effect of the accumulated information  on escape probabilities}

\subsection{General properties} Here we study properties of       $ \Bbb P\Big(\tau^b<\tau^a  \vert  X_0=x\Big)\eq  N(x;a;b)$ in terms of    parameters   $a< x<b$.  We   make the natural   assumptions on   $(X_t)_{t\ge   0} $ (see \eqref{process})
    
    \begin{ass}(A1)  Interarrival times     $\tau_n\equiv  \mathbf{t}_n- \mathbf{t}_{n-1}>0, n=1,\dots\infty$ are i.i.d. nonnegative r.v.    with  common    distribution function    $F$, namely $\tau_n\overset{\text{i.i.d}}\sim F$. 
\end{ass} 
  \begin{ass}(A2)             $(J_n)\stackrel{\text{i.i.d.}}{\sim} H$    define  an i.i.d  sequence  with a common  cdf $H$. 
 \end{ass}
 \begin{ass}(A3)   The  sequence $(J_1,\dots J_n,\dots) $ is independent of the underlying renewal process  $(N_t)_{t\ge t_0}  $.
\end{ass}

  \begin{ass}(A4)  Process   $(   N_t)_{t\ge t_0}$  has  filtration  $\mathcal N_t =\s\big(N_s, s\le t\big)$  while   $(   X_t)_{t\ge t_0}$  has  filtration $  \mathcal F_t = \s\big(X_s, s\le t\big)$. \end{ass}


\begin{proposition} Let $c>0$,  $\tau^0\eq \tau^0_-=\inf\{t>0: X_t \le   0\} $,  $\tau^b\eq \tau^b_+=\inf\{t>0: X_t \ge b\}>0 $. The function $N_b(x)$ defined by 
\beq   N: (0,b)\times\Bbb R^+\to [0,1],\  (x,b)\mapsto   N_{b} (x)= \Bbb P\Big(\tau^b<\tau^0\vert X_{0}=x\Big)\label{N} \enq 
 is monotone in both variables. If $t_b=(b-x)/c$  it  satisfies
  \begin{enumerate}  \item   \beq  \underset{ x\to b^-} \lim  N_b(x) \eq  N_{  b}(b^-)  =1\label{c>0}\enq 
  \item   $      N_{  b} (x) \le  N_{  b} (y )\le N_{  b}(b^-)  = 1, \   0< x\le y < b \label{monot} $ 
  \item   
  $ N_{  b'} (x)\ge   N_{  b} (x)\ge  N_{b\to\infty}(x)  =\Bbb P\Big(X_t>0, \forall  t\Big): =    S(x), \ x <  b'\le b\label{decr} $
 \item   
    $  \Bbb P\Big(\tau^0<\tau^b\vert X_{0}=x\Big)=1-   N_{  b} (x) \label{comp}  $
\item $    0< 1-F(t_b^-)  H(b^--x) \le  N_{  b} (x) \le  1-F(t_b^-)  H(-b)\label{bound1} $
 \item If assumptions A5,A6 below hold  $x\mapsto   N_{b} (x) $  is continuous. 
\end{enumerate}
 \end{proposition}  
 {\it Proof}.   Let $  \mathbf{U}_{ b}^x=\{ \tau^b<\tau^0\}  $ be   the event that, starting from $x$ at $t=0$,  the process escapes   $ (0,b)$   through the upper end.  Note that    $ \{\tau^{b}<\tau^{0} \}=\{X_{ {\bm \tau} }\in [b,\infty)\}$ where $ {\bm \tau}\eq \tau^0\wedge \tau^b<\infty$ w.p. $1$.    As $x$ grows so it does $(X_t)$, see  \eqref{process}  and hence the sequence $x\mapsto \{   \mathbf{U}_{ b}^x\}  $  is increasing while $b\mapsto\{   \mathbf{U}_{ b}^x\}  $   decreases as $b$ grows.  Clearly  for $c>0, x<y$
 \beq \{\tau_1 \ge t_b\}\cup \{\tau_1  < t_b, J_1\ge b-x\} \subset  \mathbf{U}_{ b}^x\subset   \{\tau_1  < t_b, J_1\le -b\}^c\text{ and }  \mathbf{U}_{ b}^x \subset  \mathbf{U}_{ b}^y \label{incl}\enq   
 which implies \eqref{bound1}. Letting  $  x\to b^- $   then  $ F(t_b^-)  H(b^--x)    \to 0$ and \eqref{c>0} follows.       

  We  next show  $\tau^\infty:=\underset{b\to \infty} \lim  \tau^b =\infty$ a.s. Indeed if    $\Bbb P(\tau^\infty <\infty)>0$  then on $\{ \tau^\infty <\infty\}$    
$$    \infty=X_{\tau^\infty} \eq x+c\tau^\infty +J_1+\dots+ J_{N_{\tau^\infty}}    \text{  a.s. which implies } N_{\tau^\infty}=\infty $$ 
namely, the  process $(N_t)$   explodes in  finite time. Besides 
  $$ \Bbb E e^{-s  (\tau_1+\dots+ \tau_{N_{\tau^\infty}}) } = \underset{n\to \infty}\lim  (\Bbb E e^{-s  \tau_1})^n=0, s>0$$    This implies   $  \tau_1+\dots+ \tau_{N_{\tau^\infty}} =\infty$ w.p. $1$ a contradiction.   Hence $ \lim \tau^\infty=\infty $    and $$    \underset{b\to\infty}\lim   \{\tau^{b}<\tau^{0} \}= \{ \tau^0=\infty \}  \text{  a.s.   }$$

  Sequential continuity of probabilities gives  (we drop  below the index $x$)
    $$ \underset{b\to\infty}\lim   N_b(x) = \underset{b\to\infty}\lim \Bbb P (  \mathbf{U}_{ b} )=  \Bbb P ( \underset{b\to\infty}\lim  \mathbf{U}_{ b})  =   \Bbb P \Big(X_t>0, \forall  t\Big)\eq S(x) $$  
 Finally since  ${\bm \tau} <\infty$ w.p. $1$. then, up to a null set, $\{\tau^0<\tau^b\}=  \mathbf{U}_{b}^c$.   (The proof of Item 6 is deferred to    Appendix A). 
\subsubsection{{\bf Symmetry properties of the escape probability}}
\begin{proposition}   \begin{enumerate}\item For $ x\in  (a,b)$, $d\in\Bbb R$ 

  \beq  \Bbb P\Big( \tau^{b+d}<\tau^{a+d} |X_0=x+d\Big)=\Bbb P\Big(\tau^{b}<\tau^{a} |X_0=x\Big) =  N_{b-a}(x-a)  \label{homog}\enq   
\item   The ``reversed'' process  $\tilde X_t=x+\tilde c t +\sum_{n=0}^{N_t}\tilde   J_n $, where $\tilde c:=-c$ and $\tilde J_n:= -J_n$  satisfies 
\beq \Bbb P\Big(\tilde \tau^{b}<\tilde \tau^{0} |X_0=x\Big) =1-N_{b}(b-x)\label{reverse}\enq 
\end{enumerate}
 \end{proposition}
 {\it Proof}.   Note  that  if there are no jumps in the time interval $(r,t)$ then 
$$  \Bbb P\Big(X_t\in dz| X_r=y\Big) =  \delta(z-y-c(t-r))dz$$
where $\delta$ is Dirac delta  with a mass at  $y+c(t-r)$. Besides
$$  \Bbb P\Big(X_{\mathbf{t}_n}\le z| X_{\mathbf{t}_n^-}=y\Big) =\Bbb P\Big( J_n \le z-y\Big)=H(z-y)$$  
Hence it is clear that $(X_t)$  is a spatially homogeneous process $-\mathbb{E}^y(g(X_t)) = \mathbb{E}^0(g(y+X_t)), \forall y-$   so, conditional on starting at $x $ the distribution of  $\tau^b$ can only  depend  on $b-x$. Choosing $d=-a$  \eqref{homog} is obtained. \eqref{reverse}  follows  noting that 
$\tilde X$  is obtained  reflecting sample paths of  $X$ over  the line $X=x$, and hence $\tilde \tau^{b} = \tau^{2x-b},    \tilde  \tau^{a} =  \tau^{2x-a} $. We finish using   item \eqref{comp}.

\begin{rem}\rm   The   invariance of  $(X_t)_{t\ge 0}$    under the  group of all space translations    and reflections permits  with no loss of generality   {\it   to suppose that 
  $a=0$ and that $c > 0$}.  
\end{rem}

\begin{proposition} Suppose $c=0$. Then condition \eqref{c>0}   {\it  needs not  to hold} as  the limit
 $  N(b)\eq \underset{x\to b^-}\lim\underset{c\to 0}\lim   N(x)\le  
\underset{c\to 0}\lim \underset{x\to b^-}\lim N(b)$ need not commute. Besides
\begin{enumerate} \item For all $x$  the bounds hold 
  $  \bar H(b^--x)\le  N(x)\le    \bar H(-x) $
 \item When    severities have a symmetric  distribution $H(x)=1-H(-x)$  then  
\beq       N(x)=   1-N(b-x) \label{symmetric} \enq  
\end{enumerate} 
\end{proposition}
 {\it Proof}.  
Noting  that 
$  \{  J_1\ge b-x\} \subset  \mathbf{U}_{ b}^x  \subset  \{  J_1 >-x\}   $ item 1 follows.    Besides if   $J_1$ has  symmetric  distribution then the law of $X$  must be invariant under reflection from the axis $z=b/2$: $  \Bbb P^x\Big(X_t\in B  \Big) = \Bbb P^{\theta(x)}\Big(X_t\in \theta\circ B \Big)   $where we  call $ \theta$  such reflection.  Since $ \theta\circ\{\tau_b<\tau_0\}= \{\tau_0
<\tau_b\}=\{\tau_b<\tau_0\}^c$ we have
$$N(x):=\Bbb P^x(  \mathbf{U}_{ b} )= \Bbb P^{\theta(x)}(  \theta\circ\mathbf{U}_{ b} )=1-  \Bbb P^{ b-x}(  \mathbf{U}_{ b} )=1-N(b-x)$$  
  \subsection{{\bf Effect of the past}}   Here we study how EP are affected by the information collected.  Let  $r $ be   a given epoch of  time (the ``present'' or `starting'' time $r$).  Clearly 
 \beq    \Bbb P\Big(\tau^b<\tau^0  \vert X_r=x\Big) =\Bbb P\Big(\tau^b<\tau^0  \vert X_0=x\Big) \label{semiM}\enq   {\it  whenever } $r\in  \Bbb T =\{\mathbf{t}_0,\mathbf{t}_1,\dots \mathbf{t}_n,\dots\}$,  the  random set of all arrival times.     However \eqref{semiM}{\it  does not extend to arbitrary present }  since     $(X_t)_{t\ge t_0}  $   needs not  being  time-homogenous nor 
 Markovian. Hence  the  escape probabilities  depend on ``starting'' time $r$ and on  which information is  
accessible.  
 In this situation  {\it  there is no real reason to   fix our attention in     $\Bbb P\Big(\tau^b<\tau^0  \vert X_0=x\Big)$    as accumulated information plays a central role}.  
   We are interested in  $\Bbb P\Big(\tau^b<\tau^0  \vert \mathcal F_r\Big)$ conditional on   the information at time $ r$,  $   \mathcal F_r=\s\Big(X_s, s\le r\Big)$.   

      Given the present $r$,       the backward and forward recurrence life   $  r\pm B_r^{\pm}   $     mark the epochs of time 
 at which the next and  last jump   occurred:   $  B_r^-=r-  \mathbf{t}_{N_r } $  and $B_r^+  =\mathbf{t}_{N_r+1}-r$.
For   $x\in\Bbb R, z\ge 0$   we introduce
\beq  N_b ( x,r,z) := \Bbb P\Big(\tau^b<\tau^0  \vert    X_r=x, B_r^-=z  \Big)\label{bigN}\enq

 \begin{proposition} {\it  For an  epoch $r$, $\Bbb P\Big(\tau^b<\tau^0  \vert \mathcal F_r\Big)$  depends only on the information contained in $ X_r$ and $ B_r^-$};  ulterior information from the past     is irrelevant.  Concretely, 

   \beq  \Bbb P\Big(\tau^b<\tau^0  \vert \mathcal F_r\Big)=  \Bbb P\Big(\tau^b<\tau^0  \vert  X_r, B_r^-\Big) =  N_b ( X_r,r,B_r^-) \label{irrelevant} \enq
 
\end{proposition}
{\it Proof}.     Note first that 
\beq \mathcal F_r=   \s\big(X_s, s\le r\big)= \s\Big( B_r^-, N_r,\mathbf{t}_0,\dots \mathbf{t}_{N_r},J_1,\dots, J_{N_r}\Big)  \enq 
\beq \text{ and } \mathcal N_r=   \s\big(N_s, s\le r\big)= \s\Big(   N_r,\mathbf{t}_0,\dots \mathbf{t}_{N_r}\Big)  \label{Nsigma}\enq  
(If $N_r=0$ we define $J_0=\mathbf{t}_0=0$).   Clearly unless  $F\sim \mathcal E(\lambda)$ neither  $(X_t)$   nor $(N_t)$ are    Markovian. Nevertheless in view of  assumptions A2,A3 and that 
\beq   X_{t+ \mathbf{t}_n} = x+ c (t+\mathbf{t}_n )  +\underset{ \mathbf{t}_j\in( 0,\mathbf{t}_n+t]}\sum J_j=X_{ \mathbf{t}_n} +ct+\underset{ \mathbf{t}_j\in( \mathbf{t}_n,\mathbf{t}_n+t]}\sum J_j  \enq 
  it follows that  given  the past of the process up to time $\mathbf{t}_n\in\Bbb T$ the future
 $(X_{t+\mathbf{t}_n}) $ is conditionally distributed as \eqref{process} starting at $ X_{ \mathbf{t}_n}$ and is independent of the past: $ \mathbb{E}^x( g(  X_{t+ \mathbf{t}_n}   ))   = 
 \mathbb{E}^{  X_{ \mathbf{t}_n}  } (g( X_t))$; 
 besides  $  Z_n:= X_{\mathbf{t}_n}, n=1,\dots\infty $  is a Markov chain.  This suggests    some  underlying simplicity.    Indeed, 
    the history previous to the last jump  is not relevant for the future evolution of  process  $X$.   At the epoch $r$ the essential history consists  only of those events of the form
  $\{X_{s }=x,  B_r^-=z, 0\le  r- z\le s\le r\}$.    More correctly, let us  define 
    $\mathcal A_r$
as  the class of   events   $$\mathcal A_r=\{X_{s_1}=x_1,\dots X_{s_n}=x_n,  B_r^-=z,  r-z\le s_1<\dots s_n\le r\}$$
for some $n\in\Bbb N$, $x_1,\dots, x_n$  and $ s_1<\dots s_n\le r$.    Then, assumptions $A1-A4$ 
  imply that 
{\it conditional on $\s(\mathcal A_r)\eq \mathcal  F_r^{\prime}$, the future evolution of $(X)$ is independent of  $\mathcal F_r$}.    
In addition,   given $X_r=x$ and $ B_r^-=z $, say, then  $X_s=x-c(r-s)$ for all $r-z\le s\le r$, i.e.  the  ``relevant'' past gets determined. (The relevant  past $\mathcal N_r$ of $(N)$ requires knowledge of both  $ B_r^- $ and $N_r$ but this does not change the argument. See \eqref{law}  below). Hence we have 

\beq \mathcal  \s(\mathcal A_r)=\s\big(  X_r,B_r^- \big)   \enq 
Note that  $\mathcal F_r$ is obtained by  joining the   sigma algebras containing  the information {\it prior and after}    the last arrival:  $\mathcal F_r = \s(\mathcal A_r) \vee \mathcal F_{r-B_r^-}$. Hence  conditional independence gives    
 \beq \Bbb P\Big(\tau^b<\tau^0  \vert  \s(\mathcal A_r) \vee \mathcal F_{r-B_r^-} \Big)=    \Bbb P\Big(\tau^b<\tau^0  \vert  \s(\mathcal A_r)\Big) =\Bbb P\Big(\tau^b<\tau^0  \vert X_r,B_r^-  \Big)  \square \enq

\section{Integral equations for the escape probability}

 \subsection{Poissonian jumps}  We consider first  the simpler case of Poisson  arrivals. 
 \begin{theorem} Suppose $c\ne 0$ and that $\tau_1\sim\mathcal E(\lambda)$. Then  $N(x)$ solves \eqref{c>0} and  
 \beq  \Big(\lambda- c\p_x \Big)N(x)= \lambda\bar H(b-x)+\lambda  \int_{-x}^{b-x}  N (x+y) dH(y),  \ -\infty<x<\infty \label{Poissongen}\enq 
\beq \text{Besides } \Bbb P\Big(\tau^b<\tau^0  \vert \mathcal F_r\Big)=  \Bbb P\Big(\tau^b<\tau^0  \vert  X_r \Big):=N(X_r) \label{Markov}\enq  
 \end{theorem}
 {\it Proof.}   Here we    take advantage that  under Poisson arrivals \eqref{process} is a L\'evy-Markov process    whose  infinitesimal  generator  $\mathbf{G} $ acts on  any $\Psi\in  L^\infty\cap C^1(\Bbb R)$   in the domain of    $\mathbf{G} $ via
\beq \Psi\mapsto \mathbf{G} \Psi (x) = c\p_x\Psi(x)+ \lambda\int_{\mathbb R} \big(\Psi (x+y) - \Psi(x)\big) dH(y)  \label{Gener} \enq  
Suppose  we allow $(X)$ to start at arbitrary $x $.   If $x\in \Bbb R-(0,b)$ then    $ {\bm\tau}\eq \tau_0\wedge\tau_b=0$ as escape occurs instantly.     Note also that $  \{ \tau^b<\tau^0\} =\{X_{ {\bm\tau} }\in [b,\infty)\}$.    This  insight yields    
\beq    
 \Bbb P^x\Big(X_{ {\bm\tau} }\in [b,\infty) \Big) =  N(x)\mathbf{1}_{(0,b)}(x)+ \mathbf{1}_{[b,\infty)}(x)+0\mathbf{1}_{(-\infty,0)}(x) \eq  \mathbf{ N} (x) \label{extendN}\enq

where $\mathbf{ N}:\Bbb R\to \Bbb R$  is an extension of  $N$ from $ (0,b )$ to $\Bbb R$. Given $ m (t,y),n  (t,y)$  let  $\Psi (t,y)\in  L^\infty\cap C^1(\Bbb R) $   solve the  Dirichlet boundary  problem 
\beq \frac{\p \Psi}{\p t} +  \mathbf{G} \Psi (t,y)=m(t,y), y\in (0,b) \text{ and }$$$$  \Psi(t,y)=n(t,y), y\notin (0,b) \label{Dirichlet} \enq 

 Additionally if   $\Psi $ is   in the domain of    $\mathbf{G} $     Dynkin's formula yields that    the  process 
\beq M_t:= \Psi (t  ,X_{t  })]- \Psi(0 ,x)- \int _{0}^{t  }  \Big(\frac{\p  }{\p s} +  \mathbf{G}\Big) \Psi (s,X_{s})\,\mathrm {d} s .   \label{Dynkin}\enq 

 is a martingale (\cite{Rolski}).  Hence       $\Bbb E   \Psi ({\bm\tau} ,X_{ {\bm\tau} })=0$ since  $ {\bm\tau}\eq \tau_0\wedge\tau_b$ is a  stopping time. Thus   for   $\Psi\in  L^\infty\cap C^1(\Bbb R)$  solving  \eqref{Dirichlet}   we have   
\beq 0=\Bbb E  M_0=\Bbb E  M_{\bm\tau}:= \Bbb E \Psi ({\bm\tau}  ,X_{\bm\tau})]- \Psi(0 ,x)-\Bbb E  \int _{0}^{{\bm\tau}  }  \Big(\frac{\p  }{\p s} +  \mathbf{G}\Big) \Psi (s,X_{s})\,\mathrm {d} s =   \label{Dynkin}\enq 
 \beq  {\displaystyle = \mathbf {E} ^{x}[n ( {\bm\tau} ,X_{ {\bm\tau} })]- \Psi(0 ,x)-\mathbf {E} ^{x}\left[\int _{0}^{ {\bm\tau} }  m(s,X_{s})\,\mathrm {d} s\right].}     \label{Dynkin}\enq

     Take   $m(t,y)=0 $ and  $ n(y)\eq \mathbf{1}_{ y\in [b,\infty)}$. 
\eqref{Dynkin} and \eqref{extendN}  give 
\beq    \Psi(x)= \mathbf {E} ^{x}[n (  X_{ {\bm\tau} })] \equiv  \Bbb P^x\Big(X_{ {\bm\tau} }\in [b,\infty) \Big)=  \mathbf{N}(x)   \enq
Thus $ \mathbf{N} (x)= \Psi(x)$   where $ \mathbf{N} $ solves 
 $   \mathbf{G}  \mathbf{N} (x) =0$ and $ \mathbf{G}$ is \eqref{Gener}.
 By insertion \eqref{Poissongen} follows $\square$ 
 
\medskip 

We consider now the general case $\tau_n\overset{iid}\sim F$. Here the {\it  above theory does not hold since  \eqref{process} is not Markov}. We resort to renewal arguments by  sharpening   the result   \eqref{incl} and ideas of section (2.2).   (Note that we drop  the $b-$ dependence in $ (\mathbf{U}_{b}^x)_{0\le x\le b}$ and  simply write  $  \mathbf{U}_{}^x  \eq  \{\tau_b<\tau_0\}$).

\begin{theorem} 
\begin{enumerate} \item      $   \mathbf{U}^x\eq  \{\tau_b<\tau_0\}$ satisfies 
\beq     \mathbf{U}_{ }^x =    \{ \tau_1\ge t_b\} \sqcup    \{  \tau_1
< t_b ,     x+c  \tau_1+J_1 \ge  b        \}  \sqcup  \Big( \{  \tau_1< t_b ,   0\le x+c \tau_1 +J_1 < b        \} \cap   \tilde{\mathbf{U}}^{z}   \Big) \label{decom2} \enq
i..e it   can be decomposed as a disjoint   union  where     $ \tilde{\mathbf{U}}^z $ satisfies  for given $z$ 
\beq   \tilde{\mathbf{U}}^z      \perp\kern-6pt\perp      \s\Big(  \tau_1,J_1\Big)
 \text{ and } \   \tilde{\mathbf{U}}^z\overset{\Bbb P}=
\mathbf{U}_{ }^z \text{  or } \Bbb P( \tilde{\mathbf{U}}^z) = \Bbb P(\mathbf{U}_{ }^z), \forall z\label{Aprop}\enq 

  \item Let    $   \mathbf{U}_{ r}^x $ be  the event that $(X)$   exits through the upper barrier when      $X_r=x$ and $r>0$ is the present time.  Call  $ J\eq J_{N_r+1}$.   Then $  \mathbf{U}_{ r}^x$   can be decomposed as the disjoint   union 
 $        \mathbf{U}_{ r}^x =     \mathbf{U}_{ r}^{x,(1)} \sqcup    \mathbf{U}_{ r}^{x,(2)}  \sqcup  \mathbf{U}_{ r}^{x,(3)} \text{ where  }  $
$$  \mathbf{U}_{ r}^{x,(1)} \eq \{ B_r^+\ge t_b\},   \  \mathbf{U}_{ r}^{x,(2)} \eq \{   B_r^+
< t_b ,   0 < x+cB_r^+ +J  < b        \}   $$  
\beq    \text{ and } \mathbf{U}_{ r}^{x,(3)} \eq \{   B_r^+
< t_b ,   0 < x+cB_r^+ +J  < b        \} \cap      \tilde{\mathbf{U}}  \label{decom} \enq  is made up of independent events.  

If $r+B_r^+\eq  \mathbf{t}\in\Bbb T$ then  $   \tilde{\mathbf{U}} = \mathbf{U}_{ \mathbf{t} }^{    X_{  \mathbf{t} } } $ is conditionally  independent of $\mathcal F_{\mathbf{t}}$ given  $X_{  \mathbf{t} } $ and  \beq \Bbb P(\mathbf{U}_{ \mathbf{t} }^{    X_{  \mathbf{t} } }|    \mathcal F_{\mathbf{t}})=\Bbb P(\mathbf{U}^{ z } )|_{z= X_{  \mathbf{t} }   }\label{condind} \enq  \end{enumerate}
\end{theorem}

 {\it Proof}. We prove \eqref{decom} since then \eqref{decom2} follows  letting $r=0$ and noting that 
 $B_0^+=\tau_1$.  After time $r$, given     $X_r=x$,     four excluding possibilities  unfold, depending on the  evolution   up to the first arrival: $s\mapsto X_s, r\le s\le  r+ B_r^+ $  
 \begin{enumerate}\item 
   $ B_r^+\ge t_b $. Then $ X_{r+t_b}=b$.

If, by contrast,   $ B_r^+ < t_b$ then  a jump $J$ occurs at $r+B_r^+$ prior to escape. Then \item  
  $ B_r^+ < t_b, J>0 $   and  
$     X_{  r+ B_r^+}\eq x+c B_r^++J  \ge  b $.  Then  $ X_{  r+ B_r^+} \ge  b$.

 \item  $ B_r^+ < t_b $  and  $   0\le  x+c B_r^++J < b$. After the ``first'' renewal $r+ B_r^+ $  the process starts  at  $ X_{  r+ B_r^+}\eq x'\in (0,b)$ and will   exit through the upper barrier if      $   \mathbf{U}_{ \mathbf{t} }^{    X_{  \mathbf{t} } } $ occurs.

In all these cases escape  will occur  through the upper  barrier.  
  \item  
  $ B_r^+< t_b, J<0$   and  $  X_{  r+ B_r^+} <0$.   Then  $(X)$
 escapes  through the lower barrier. \end{enumerate}

   This implies \eqref{decom} where  $ \tilde {\mathbf{U}}^{  }  \equiv \mathbf{U}_{ \mathbf{t} }^{    X_{  \mathbf{t} } } $ where    $\mathbf{t}\eq r+B_r^+\in\Bbb T, z\eq  X_{  \mathbf{t} }$.  
We now see \eqref{decom} and \eqref{condind}. Let   $z=  X_{ \mathbf{t} },   t'_b=( b-z)/c$. Then   
$$   \mathbf{U}_{ \mathbf{t}}^{ z} =      \{  \tau_{N_r+2}\ge t'_b\} \sqcup    \{  \tau_{N_r+2}
 < t'_b ,     y+c   \tau_{N_r+2}+J_{N_r+2} \ge  b        \}  \sqcup \dots \in \s \Big(\tau_n, J_n, n\ge N_r+2   \Big) $$
By contrast $$  \{ B_r^+
 < t_b ,   0 <x+cB_r^+ +J < b        \} \in   \s \Big(\tau_{N_r+1}, J_{N_r+1} \Big)   $$
The result follows since assumption A3 implies $\s \Big(\tau_{N_r+1}, J_{N_r+1} \Big)\perp\kern-6pt\perp  \s \Big(\tau_n, J_n, n\ge N_r+2   \Big)  $ where   $\mathcal G\perp\kern-6pt\perp  \mathcal F$  denotes independence of $\s-$ fields.    
Note also that  $\tau_{N_r+2}  \overset{F}=\tau_1, $     $J_{N_r+2}\overset{H}=J_1,\dots $ have the same law   (see A1, A2); hence  $$   \mathbf{U}_{ \mathbf{t}}^{ z} \overset{\Bbb P}=      \{ \tau_1\ge t'_b\} \sqcup    \{  \tau_1
 < t'_b ,     y+c  \tau_2+J_1 \ge  b        \}  \sqcup\dots = \mathbf{U}_{}^z     \quad \square$$

  Actually   
   $\Bbb P\Big(\tau^b<\tau^0  \vert   \mathcal F_r \Big)$,  given by 
 \eqref{irrelevant},   is   retrieved once  $N_b(x)$ is known.  

In the sequel we make the mild and convenient assumptions
\begin{ass}(A5)   {\it  $t_b$ is not a   mass of }  $F$.
\end{ass}  
  \begin{ass}(A6)    $ L:= \Bbb P\Big(\tau_1\le b/c, J_1\in (-b,b)\Big)<1\label{cont} $
\end{ass} 
 \begin{ass}(A7)$F$ has a density $f$ and  $ H$ has a density $h$.\end{ass} 
Thanks to  A5 we avoid the messy distinction between  $F(t_b^-)$ and $F(t_b)$  while A6 guarantees that the IE \eqref{Fred} below  satisfies a fixed point condition. Assumption A7 is unnecessary  at this stage, but will be convenient  when we take up the task of solving \eqref{Fred}  (sections 4-6). 

We start considering the conditional distributions of the    Markov process  $t\mapsto B_t^+$.  
 \begin{Lemma} For any epoch $r$, $ B_r^+ $ is   conditionally independent of $X_r$ and the history $\mathcal F_r $  given $ B_r^- $. Besides   
\beq  {\displaystyle  (B_r^+ \perp \!\!\!\perp  \mathcal F_r )\mid B_r^-  }, \quad \text{   }  B_r^+  \perp\kern-6pt\perp  J_{N_r+1} \text{ and }  J_{N_r+1} \perp\kern-6pt\perp\mathcal F_{N_r}  \enq
\beq \Bbb P \Big( B_r^+ > t\vert  \mathcal F_r\Big)=\Bbb P \Big( B_r^+ > t\vert  \mathcal N_r\Big)=  \Bbb P \Big( B_r^+ > t\vert B_r^- \Big)=  \bar F( t+z) /\bar F(z)|_{z= B_r^-}\label{law}\enq 
     
\end{Lemma}
{\it Proof}.      Given $ B_r^-  $  assumption A1 with random  $n(\omega)\eq N_r(\omega)$  yield  
$$B_r^+ =  \tau_{N_r+1} -  B_r^-\perp\kern-6pt\perp    \mathcal N_r=\s(N_r,\tau_0,\dots \tau_{N_r})$$     Thus  $ B_r^+ $ is      conditionally independent of the history $\mathcal N_r $  given $ B_r^- $.     $ (X)$  contains  information on $B_r^-$;  hence  $B_r^+ ,  X_r$ are not independent, but they are given $B_r^-$. It follows that  
 $$  \Bbb P \Big(  \tau_{N_r+1}  > t+z\vert B_r^- =z\Big)=   \sum_{n=0}^\infty \Bbb P \Big(  \tau_{N_r+1}  > t+z\vert B_r^- =z, N_r=n\Big)   \Bbb P \Big(    N_r=n\vert B_r^- =z\Big)  $$ $$ \sum_{n=0}^\infty \Bbb P \Big(  \tau_{n+1}  > t+z\vert     \tau_{n+1}>z \Big)   \Bbb P \Big(    N_r=n\vert B_r^- =z\Big)  =\frac{\bar F( t+z) }{\bar F(z)}  \sum_{n=0}^\infty   \Bbb P \Big(    N_r=n\vert B_r^- =z\Big) $$  
 $$ \text{ or }  \Bbb P \Big( B_r^+ > t\vert B_r^- =z\Big)= \bar F( t+z) /\bar F(z) $$  

\begin{theorem}   \label{ineq}   $N_b(x,r,z)$ and $\Bbb P\Big(\tau^b<\tau^0  \vert   \mathcal F_r \Big)$   follow  from  $N_b(x)$ via

  $$N_b(x,r,z) =  \frac{  1}{\bar F(z)}\Big(     \bar F (z+t_b) + \int_0^{t_b}\bar   H(b-x-cl )  d_lF(z+l) +$$  \beq  
  \int_{ 0}^{ t_b}     d_lF(z+l)  \   \int_{-x-cl}^{b-x-cl}  dH(y)      N_b(x+cl+y) \Big)\quad  \label{hist} \enq

In particular $N_b ( X_r,r,B_r^-) \eq N_b ( X_r,B_r^-)  $ is independent of the present epoch $r$.

\end{theorem}



        In the sequel to ease  notation  we set   $J\eq J_{N_r+1}$ while  
$\vec X=\vec x $ stands for $ X_r=x,  B_r^-=z $. It follows from \eqref{irrelevant} and 
Proposition (4)  that $$N_b(x,r,z)=   \Bbb E\Big(   \mathbf{1}_{ \tau^b<\tau^0} \vert X_r=x,B_r^-=z\Big) =\sum_{j=1}^3 N_b^j(x,r,z)$$ 
where   $ N_b^j(x,r,z)$   denote the probabilities   of the different terms  appearing in the RHS of   \eqref{decom}.
Hence  \eqref{law} gives $  N_b^1(x,r,z)=  \Bbb P \Big( B_r^+>t_b\vert B_r^-=z\Big) = \bar F(z+t_b)/\bar F(z)$.

        We evaluate next the probability of $\mathbf{U}_{ r}^{x,(3)}$  using the tower property as
 $$ N_b^3(x,r,z)\eq \Bbb P\Big(    \mathbf{U}_{ r}^{x,(3)}  \big|  \vec X=\vec x \Big)=\Bbb E\Big( \Bbb P\Big(    \mathbf{U}_{ r}^{x,(3)}      \big|   J , B_r^+,   \vec X=\vec x  \Big)\vert   \vec X=\vec x
\Big) =$$$$    \int  \Bbb P\Big( \mathbf{U}_{ r}^{x,(3)} \big|   J=y,B_r^+=l ,   \vec X=\vec x \Big)  \Bbb P\Big(  J\in  dy,B_r^+\in dl   \big| \vec X=\vec x\Big)=$$ 
     $$  
    \int \Bbb P\Big(    \mathbf{1}_{   l \le  t_b }   \mathbf{1}_{\{  -x-c   l\le  y <b- x-c l      \}} \cap   \mathbf{U}_{  r+B_r^+}^{ X_{  r+ B_r^+} }  \big|   J=y, B_r^+=l ,   \vec X=\vec x \Big) \cdot$$$$  \Bbb P\Big(  J\in  dy,B_r^+\in dl   \big|  \vec X=\vec x\Big)$$ 
 Besides  if $  B_r^+= l \le  t_b, J=y$ then $X_{  r+ B_r^+}=x+cl+y$. Hence \eqref{condind} gives  
   $$ \Bbb P\Big(       \mathbf{U}_{  r+B_r^+}^{ X_{  r+ B_r^+} }    \big|   J=y, B_r^+=l ,    \vec X=\vec x \Big) =$$
 \beq\Bbb P\Big(     \mathbf{U}_{\mathbf{t}_{N_r+1}}^{x+cl+y }  \big|   J=y, B_r^+=l ,   \vec X=\vec x \Big)  =\Bbb P (\mathbf{U}^{  x+cl+y}_{ 0}   ) =N(x+cl+y) \label{NN}\enq 
    y  conditional  independence (Lemma 1)  and  assumptions  A1-A3
      we have 

 $$ \Bbb P\Big(  J_{N_r+1}\in  dy,B_r^+\in dl   \big| \vec X=\vec x  \Big)= \Bbb P\Big(  J_{N_r+1}\in  dy  \big| \vec X=\vec x  \Big) \Bbb P\Big(  B_r^+\in dl   \big|   X_r=x,  B_r^-=z \Big)   $$ \beq  =\Bbb P\Big(  J_1\in  d y\Big)   \Bbb P\Big( B_r^+\in dl   \big|  B_r^-=z \Big)  \label{Fin}\enq  
Hence eqs. \eqref{law}, \eqref{NN} and \eqref{Fin} give 
$$ \Bbb P\Big(    \mathbf{U}_{ r}^{x,(3)}  \big|  \vec X=\vec x \Big)=   \int_{ 0}^{ t_b}     \frac{d_l F(z+l)}{\bar F(z)}  \   \int_{-x-cl}^{b-x-cl}  dH(y)      N_b(x+cl+y)     \quad \square $$

 \subsection{{\bf Integral equations, general case}}\label{decomp}We have seen how to codify EP conditional in the history in terms of the basic object   $N_b(x)$. Here \eqref{Poissongen} is not valid. 
 We now derive integral equations for  this object. 
\begin{theorem}  Suppose assumptions A1-A4 hold. 
  The function $ x\mapsto N_b(x)\eq N(x)    $ of \eqref{N}   satisfies for $ 0\le x\le b$  the IE 
\vskip0.2cm
  \beq N(x)= \bar F (t_b)+    \int_0^{t_b} \Big(     \bar   H(b-x-cl )+ 
    \   \int_{-x-cl}^{b-x-cl}  dH(y)      N(x+cl+y)\Big)dF(l)   \quad  \label{Fred}  \enq  

\end{theorem} 
{\it Proof}.  If $r=0$ is  $z=B_0^-=0,  B_0^+=\tau_1, N_b(x,0,0)= N_b(x) \text{ and } $  $ \Bbb P\Big(  B_0^+
  > l   \big| B_0^-=0  \Big)= \Bbb P (  \tau_1> l  )=\bar F(l)$.  
 Inserting these values in    \eqref{hist}  we  find \eqref{Fred}. 

 \begin{Corollary} 
 \eqref{Poissongen} and  \eqref{Markov}  hold if  either  $\tau_1$ is   exponential        or  
   $c=0$.  If  
    $r\in\Bbb T$ is one of the jump times then \eqref{Markov} holds.  
\end{Corollary}   
{\it Proof}.  
 If $c=0$ then $t_b=\infty$ and  both \eqref{Fred} and  \eqref{hist}  simplify to 
\beq  N_b(x,r,z) =          \bar   H(b-x )   + 
    \int_{-x}^{b-x}       N_b(x+y) dH(y) =N_b(x)\label{26} \enq   

Finally,  if $r\in\Bbb T$ $B_r^-=0 $ and $ \Bbb P\Big(  B_r^+
  > l   \big| B_r^-=0  \Big)= \bar F(l) \quad \square$
 \vskip0.2cm


         To clarify the role of different contributions appearing in  \eqref{Fred}  we  
 decompose $H$ a convex  combination   of {\it proper } distribution functions $H_{j,\pm}$: 
 \beq H = q_1 H_{1-} + q_2 H_{2-}  + p_1 H_{1+} + p_2 H_{2+} \label{H}   \enq
     $$\text { where } p_j,q_j\ge 0,  p_1+p_2+q_1+q_2=1 \text{ and } $$   
      $$\text{   supp }  H_{2-}= (-\infty, -b], \text{ supp }  H_{1-}= (-b,0)  \text{ supp }  H_{1+} = [0, b-x ),\text{and  supp }  H_{2+}=[b-x,\infty)$$ 
\begin{proposition}\label{trivial}\begin{enumerate} 
 
 \item   
  If      support  $H=(-\infty,-b]\cup [b-x,\infty)$, i.e.  $p_1=q_1=0$,  then 
  \beq N_b(x)= \bar F (t_b)+  p_2 F (t_b)  \quad  \label{double} \enq  
\item Let  $N_b(x)$, $\tilde N_b(x)$ be  the  EP  corresponding to  jump cdf's $H$ and, respectively,  $\tilde H$. If     $\tilde H=H$  on $ ( -b, b-x )$ and $p_2=\tilde p_2$  then  $N_b(x) =\tilde N_b(x)$.  
 
\end{enumerate} \end{proposition} 
   

 {\it Proof}. Using that     $ H(y)= q_1+q_2+p_1 H_{1+}(y ) \text{  if  }     0\le y <b-x, $

 $ \    H(y)=q_1 H_{1-}(y)+ q_2       \text{  if  } -b < y\le 0   \text{ and that }   $  $$0\le l<t_b\Rightarrow 0\le b-x-cl\eq y<b-x, \bar H(y)= p_2+ p_1  \bar     H_{1+} (y)$$   \eqref{Fred} reads  
 $$ N(x)=   p_1 \int_0^{t_b}  dF(l)  \Big(  \bar     H_{1+} (b-x-cl )+ 
   \int_{0}^{b-x-cl}     N(x+cl+y) d       H_{1+}(y)\Big) $$ 
\beq  + q_1  \int_{ 0}^{ t_b}   dF(l) \int_{-x-cl}^0  N(x+cl+y) dH_{1-}(y)+\bar F (t_b) + p_2  F (t_b)\label{fav2}  \enq 
  
    Setting   $p_1=q_1=0$ \eqref{double}  follows. Alternatively, note that  when $p_1=q_1=0$ \eqref{decom2}  reads  $ \{ \tau^b<\tau_0\}=  \{ \tau_1\ge t_b
\}\cup  \{ \tau_1 < t_b , J_1>0 
\}     $.   \eqref{double}  follows trivially. 

    The structure of \eqref{fav2} shows that  that   if    $\tilde H_{1\pm}= H_{1\pm} $ and $ p_2=\tilde p_2 $  then  $N$ and $ \tilde N$ solve the same equation and hence $N= \tilde N$. 

\begin{rem}  \rm In Appendix A we prove that  under the mild assumption  A6      \eqref{Fred}  involves a     Lipschitz continuous   operator   with Lipschitz constant $  L<1$   and hence  has a unique solution.  Hence {\it  the problem of obtaining the EP is codified in a one-to-one way into solving  such   linear Fredholm  IE}.   Since   generically     {\it  a closed form solution    is not possible} to help a reader place properly   the  situation   we indicate   the decomposition \eqref{H}  of $H$ as  $\vec p=( q_1, q_2,p_1,p_2)$.   We study   assumptions that render \eqref{Fred} solvable in terms of the factors present in decomposition \eqref{H}, denoted $\vec p\eq (q_1,q_2,p_1,p_2)$.   Case $p_1=q_1=0$ is trivial, cf. \eqref{double}.  \end{rem}

 \section{ Actuarial case:   $\vec p=(\star, \star,0,0)$} 
 Here we discuss the situation of the  risk   model  when   {\it only   negative jumps are allowed}. This corresponds to $p_1=p_2=0$, or  $\vec p=(\star, \star,0,0)$ where $\star$ marks  non-null components. We shall   determine a  {\it   large class of arrivals distributions  under which this equation is solvable}.  Concretely   A7 is  modified to
 \begin{ass}(A8) Supp. $J_1 \subset(-\infty,0) $. Besides     $\tau_1$ and $ -J_1$ have densities $f_{\tau_1}(t) \eq f(t)$ and $f_{-J_1}\eq h$.  \end{ass}     
We can write \eqref{Fred} as
\beq N(x)=   \bar F (t_b) +  (1/c) \int_{x}^b   dl   f((l-x)/c ) \int_{0}^{ l} N( l-y) h(y) dy \label{Fredac2} \enq 




\subsection{{\bf Poisson arrivals}} 
 Here we consider    the risk   model with    Poisson  arrivals  $f(t) = \lambda   e^{-\lambda t} $.     Eq. \eqref{Poissongen}  yiedls that $N$ satisfies also the   integro-differential equation
\beq\Big(\lambda-c\p_x \Big) N(x)= \lambda \int_{0}^{ x} N( x-y) h(y) dy, \ 0\le x\le b\ \label{N=1}\enq

The key observation  is that|unlike \eqref{Fredac2}|  \eqref{N=1}  has constant coefficients and is well suited to  Laplace transformation by appropriate  extension to    the entire line $0\le x <\infty$.   If  $ \Upsilon(x)$ denotes a {\it general } solution  to \eqref{N=1} its Laplace transform (LT) must satisfy 
$$\Big(  (\hat h(s)-1)\rho+s \Big) \hat \Upsilon(s)= \Upsilon_0, \   s\ge 0 $$

where   $\rho=\lambda/c$ and $ \Upsilon_0 \eq \Upsilon(0) $. 
  By the Laplace inversion formula $ \Upsilon(x) $ can be recovered as $ \Upsilon(x)\equiv \Upsilon_0 \pi(x)-$where  \beq \pi (x)\eq \frac{1}{2\pi i}\int_{\beta-i\infty}^{\beta+i\infty}\frac{e^{sx}ds}{(\hat h(s)-1)\rho+s}\label{escapeexp} \enq
 and the line $s=i\beta, \beta >0$ lies to the right of all singularities, namely we take  the standard Bromwich contour. 
 We prove below (see \eqref{norm}) that  $\pi(0)=1$ and hence the consistency condition: $ \Upsilon_0 \eq \Upsilon(x=0) $ is identically satisfied. Hence $ \Upsilon(x)$ defines a one-parameter family of solutions labeled by the {\it free constant} $ \Upsilon_0 $. To  retrieve $N_b(x) $  we must require    the extra boundary condition  \eqref{c>0}. This gives:

\begin{theorem} Suppose A8 holds,  $\tau_1\sim\mathcal E(\lambda)$ and  $\pi(x)$ is given by \eqref{escapeexp}. Then the  EP    is
\beq N_b(x) \eq \Bbb P^x\Big(\tau^b<\tau^0\Big) = \pi(x)/\pi(b) \text{  } \label{escape2}\quad \square  \enq

\end{theorem}

We now elaborate further on the meaning of  $\pi(x)$  and  recover the ruin probability. 
\begin{proposition}\begin{enumerate} \item $\pi(0)=1$
\item  Suppose  $h(x):= \Phi(x-b)\theta(x-b)$ vanishes for $0\le x <b$.    Then  $  \pi(x)=    
e^{\rho x }         $ and (see \eqref{double})  $  N(x)=    
e^{\rho (x-b) }        $. Here and elsewhere $\theta(x)=\mathbf{1}_{x\ge 0}$ is the Heaviside function.
 \item Let  $ m:=-\Bbb EJ_1$. If    $\rho m\ge 1$     the survival  probability vanishes: $ S(x):=  N_{b\to\infty}(x) =0$.  

  If $\rho m<1$  then 
  $\underset{x\to\infty}\lim \pi(x)=  (1- \rho m )^{-1} $ and  $S(x)$   satisfies
 \beq  S(x) =   \pi(x)/\pi(b\to\infty)=    \frac{1-\rho m}{2\pi i}\int_{\beta-i\infty}^{\beta+i\infty} 
      \frac{e^{sx}ds}{(\hat h(s)-1)\rho+s}  \label{ruin2}\enq 
 The EP satisfies also 
\beq    \Bbb P^x\Big(\tau^b<\tau^0 \Big)  =  \frac{ S(x) } { S(b) }\label{escape3} \enq

\end{enumerate}  \end{proposition} 

{\it Proof}. We can evaluate $\pi(0)$ as \beq  \pi(0)= \frac{ 1}{2\pi i}\int_{\beta-i\infty}^{\beta+i\infty} \frac{ds }{s}
\Big(1- \frac{ \rho( \hat h(s)-1) } { s+ \rho( \hat h(s)-1) } \big)=1 \label{norm}\enq


 Indeed, since $\hat h(s)$ is analytic on $\Bbb C^+:=\{s=s_R+is_I\in\Bbb C: s_R\ge 0\} $,    the second   integrand is $O(1/s^{ 2})$ and is analytic  on     $s_R \ge 0$. We  close the contour with  a very large semi-circle   in the  right half-plane whereupon it vanishes by    Cauchy's theorem. 

   For item (2) note that   $\hat h(s)= \hat \Phi(s) e^{-bs}$  and hence   $$ 2\pi i\pi(x)=  \int_{\beta- i\infty}^{\beta+i\infty}     \frac{  e^{sx} ds} { s+\rho  \hat \Phi(s) e^{ -b s}  -\rho}= \sum_{n=0}^{\infty}\sum_{j=0}^{n} \int_{\beta- i\infty}^{\beta+i\infty} \frac{ds}{s}      e^{sx}  (\rho/s)^n \hat \Phi^j(s) \binom{n}{j}  e^{ s(x-bj)} $$ 
  
  We close
 Bromwich  contour  on the  on the right   half-plane   if $j\ge 1-$  whereupon the integral vanishes by analicity.     If  $j=0$  we close on the left. The residue theorem gives  $\int_{\beta- i\infty}^{\beta+i\infty}        ds (2\pi i s)^{-1}  s^{k-n}=\delta_{nk}$ and hence
 $$ \pi(x)=    \sum_{n=0}^{\infty} \int_{\beta- i\infty}^{\beta+i\infty}       ds (\rho/s)^n   (2\pi i s)^{-1}  \sum_{k=0}^{\infty}    (sx)^k/k!  = \sum_{n=0}^{\infty}  (\rho x)^n/n! =  e^{\rho x} $$  
 
  For  (3) note that     the following conditions are equivalent C1: $\pi(x)$ is   bounded.   C2: $S(x)>0$ (see \eqref{escape2}).   C3:    The net profit  condition (NPC) holds (see \cite{Mikosch}): \beq \Bbb E \tau_1.|\Bbb E J_1|<\infty, \  c\Bbb E \tau_1+\Bbb E J_1 > 0 \label{NPC}\enq 
 C4:         All non-null zeroes of    Lundberg equation (LE)  $s=   (1-\hat h(s))\rho$  have negative real parts (a complete discussion of  this aspect is performed in   \cite{Berger2,lg04}).      Hence if   $\rho m<1$   we   can appeal to dominated convergence  and   the  Tauberian final value theorem for Laplace transforms (\cite{Feller3}) to find    $\underset{s\to 0}\lim    s (\hat h(s)-1) =m\eq -\Bbb E J_1>0 $ and  
  \beq \underset{x\to\infty}\lim \pi(x)=\underset{s\to 0}\lim    s\hat \pi(s)=\underset{s\to 0}\lim    s\Big((\hat h(s)-1)\rho+s\Big)^{-1}= (1- \rho m )^{-1}, \label{fin}\enq

     The  survival  probability follows  recalling $ S(x):=  N_{b\to\infty}(x)$.  
 
\begin{rem}\rm     
While result   \eqref{ruin2} for  {\it survival  probabilities} is well known  the connection with escape probabilities   \eqref{escape2} is a different matter.  In a seminal paper  Bertoin  \cite{Bertoin} introduces a  similar  factorization   for stable one-sided L\'evy process  in terms of the L\'evy measure. This line of approach  is continued in \cite{Surya,Chiu,Avram}.  We note that         representation \eqref{escape3} does not hold if  either   assumption    $\tau_1\sim\mathcal E(\lambda)$ or $cJ_1<0 $ is  dropped  as we show    below.  
Note also that {\it the scale function  of a difussion }  (cf.  \cite{Feller1,Feller3}) 
 is any strictly increasing  $s :  (a,b)\to \Bbb R$   such that
\beq  \Bbb P^x\Big(\tau^b<\tau^a\Big) =\big(s(x)-s(a)\big) /\big(s(b)-s(a)\big)\label{scale}\enq    

 \end{rem}  
 
\begin{rem}\rm    An analogue  of  the   Pollaczek-Khinchine   formula  for   escape  probabilities follows by   
   series  expansion in the Laplace integral \eqref{escapeexp} 
$$\Big((\hat h(s)-1)\rho+s\Big)^{-1} 
=   s^{-1}  \sum_{n=0}^\infty (\rho  \hat {\bar H})^n$$   \end{rem}

 \subsubsection{Several examples}\begin{exmp}Exponential jumps \rm    $-J_1\sim \mathcal E(\gamma)$. Then $N_b(x)= \frac{  1+ \gamma      x   }{   1+ \gamma      b  }   $ when $ \rho =\gamma   $ and  
\beq   N_b(x) =    \frac{  1 -(\rho/\gamma)  \e^{(\rho-\gamma)x}}{  1 -(\rho/\gamma)   \e^{(\rho-\gamma)b}}, \  \rho
  \ne \gamma \label{Poissexp}  \enq    


 \end{exmp} \begin{exmp} Rational jumps. \rm  Suppose that  $h(s)$  has rational LT $ \hat h (s)$.
Then  $\pi(x) $ is given by sums of terms of the form $\pi(x)\approx \sum_i  \gamma_i e^{s_i x} $ where $s_i\in \mathcal R$ are the zeroes of the denominator.  
 \end{exmp}
\begin{exmp}Gamma jumps.  \rm When  $-J_1\sim \Gamma(\alpha,\gamma), \alpha\notin \Bbb N$  then     $ \hat h (s)=   \Big(\frac{  \gamma  }{   \gamma+  s  }\Big)^\alpha   -$ and hence  $ \hat \pi (s)-$   are  non-rational  and have a branch cut  singularity on $(-\gamma,0)$. Nevertheless in particular cases it is   possible to perform the inversion of \eqref{escapeexp}. Here we suppose  $ \alpha=1/2 $.  We  rationalize  $ \hat \pi (s)$   to the convenient form 
        \beq  
          \hat \pi(s) = \frac{ \sqrt {   \gamma+  s  } }{  (\sqrt{  \gamma  }  -\sqrt {   \gamma+  s  } )\rho+s\sqrt {   \gamma+  s  }} =\frac{(s+\gamma) (s-\rho) -\rho \sqrt { \gamma }\sqrt {s+\gamma }    } { s(s^2 +(\gamma-2 \rho)s- \rho(\rho-2\gamma)  }. \label{IG1/2Laplace} \enq
Let $2s_{\pm}=2\rho- \gamma\pm \sqrt{\gamma(\gamma+4\rho)}$ and $\xi_{\pm}\eq s_{\pm}+\gamma$. Using a partial fraction expansion lengthy calculations yield $\pi(x) =   S(x)/(1-\rho m)$ where $S$ is given in terms of  erf function as
 \beq 
S(x)= \frac{1}{2}\Big(1+ \text{ erf } (\sqrt{\gamma x} )\Big)  + \frac{  1}{2\gamma  \rho (s_--s_+)}  \Big(e^{s_+x}s_- \xi_+(s_+-\rho) - e^{s_-x}s_+ \xi_-(s_--\rho)  \Big) +
 $$
$$\frac{ 1}{2  \sqrt{\gamma} ( s_--s_+)   }      \Big(\sqrt{\xi_-}s_+ e^{s_-x} \text{ erf } (\sqrt{x\xi_-}) - \sqrt{\xi_+} s_- e^{s_+x}  \text{ erf } (\sqrt{x \xi_+})
    \Big) \label{Gamma}\enq   \end{exmp}

 \begin{exmp}  Constant jumps $-J_1=y_1\in\Bbb R^+$.\rm  This   corresponds to an actuarial situation   where policy holders have the right to a {\it  fixed predetermined } compensation $y_1 >0$ per claim.  Hence  $\hat h(s)=e^{  -y_1 s} $ and 
  
   $$  \pi(x)= 1/(2\pi i) \int_{\beta- i\infty}^{\beta+i\infty}     \frac{  e^{sx} ds} { s+\rho e^{ -y_1 s}  -\rho} =C(x,-\rho,1) \text{ where} $$  
   $$C(x,\rho,p) =1/(2\pi i) \int_{\beta- i\infty}^{\beta+i\infty}    \Big( \rho+s- p\rho  e^{- y_1 s}    \Big)^{-1} e^{sx }   ds =     $$ 
$$=    \int_{\beta- i\infty}^{\beta+i\infty}   \frac{ds }{2\pi i s}  \sum_{k=0}^\infty   \sum_{n=k}^\infty  ( -  \rho /s)^n \binom{n}{k}     e^{ s(y-k y_1)}   (-p)^k   $$
The integral is evaluated  by series expansion, where term-wise we close
Bromwich  contour with a large half circle on the left or right half $s-$plane. Only when   $   x-k y_1\ge 0$  we can  close on the
right half-plane and pick a pole at $ s = 0$. Call   $k_1\eq  \lfloor     x/y_1 \rfloor $.  The residue follows from 
 $$=    \int_{\beta- i\infty}^{\beta+i\infty}\frac{ds }{2\pi i s}  \sum_{k=0}^{k_1}   \sum_{n=k}^\infty   \sum_{j=0}^\infty    ( -  \rho /s)^n \binom{n}{k}        (-p)^k    s^j(x-k y_1)^j/j!  $$ 
 $$=    \sum_{k=0}^{k_1}     \sum_{n=k}^\infty      ( - \rho  )^n      (-p)^k       (x-k y_1)^n/((n-k)!k!)  =     
      \sum_{k=0}^{k_1}    ( p\rho  (x-k y_1))^k  e^{  -   \rho (x-k y_1)}          /   k !  $$
  
 Hence $\pi(x) =C(x,-\rho,1)= e^{\rho x } \tilde \pi(x)  \text{ where }$
\beq    \tilde \pi(x)\eq   \sum_{k=0}^{k_1}  a_k(x) \text{ and  } a_k(x)=  (-\rho e^{ -y    \rho }    )^k   (x-ky_1)^{k } /k!  \label{cost-}   \enq  
As  $x\to\infty$ is        $k_1\to\infty$ and convergence of the series is unclear. Note first that  $ \tilde \pi(x)$   
 involves an {\it  alternating series} whose general term $a_k\downarrow 0$ monotonically as $ k\to\infty$. Thus $\tilde \pi(x)$ converges as  $x\to\infty$. The convergence of $\pi(x\to\infty)$ is delicate: the ratio test  shows that 
  it  requires   $y_1    \rho e^{ 1-y_1    \rho }<1\Leftrightarrow 1-\rho y_1>0$.  In this case NPC and  \eqref{fin} hold. One has $ \underset{x\to\infty}\lim  \pi(x)  = (1-y_1\rho)^{-1} $

  \end{exmp}

   \subsection{{\bf Erlang $ \Gamma(n,\lambda)$ and hypo-exponential arrivals}}   Here we generalize the previous results   {\it   to    the actuarial model  under    hypo-exponential arrivals. } That is, A8 holds and  there exist      parameters  $0<\lambda_1\le \dots \lambda_n$   such that  $f_{\tau_1}  \eq f$ satisfies
 \beq \hat f(s)= \prod_{  j=1}^n \frac{ \lambda_j }{(\lambda_j+s) }\label{hypoexp}\enq
   \begin{rem}\rm This    distribution 
       corresponds to a   sum   $\tau_1=X_1+\dots X_n$ of $n$ independent variables $X_i\sim\mathcal E(\lambda_i)$.
      Interesting particular  cases are \begin{enumerate} \item $ \lambda_1= \dots =\lambda_n$: this yields    Erlang distribution $ f_{\tau_1}(t)= \lambda_1 ( \lambda_1 t)^{n-1} e^{-\lambda_1 t}/(n-1)!$.  
\item   $\lambda_j=j\lambda_1,\forall j$. Here   $\Bbb{P}(\tau_1 \leq t) = 
  (1-e^{-\lambda_1 t} )^n $, the  {\it order statistics sampled from an exponential distribution.}
\item     Under    strict generic   inequalities   the density is the Lagrange combination 
  \beq   f(t)=\sum _{k=1}^{n}p_k \lambda _{k}e^{-t\lambda _{k}}
\text{ where } p_k\eq  \prod _{j=1,j\neq k}^{k}{\frac {\lambda _{j}}{\lambda _{j}-\lambda _{k}}}  \in\Bbb R \enq
 \end{enumerate}\end{rem}   This situation modeled by   \eqref{Fredac2}  is not solvable as stands; nevertheless {\it   one can transform it to an equivalent, but simpler  integro-differential equation    with  appropriate   BCs} that generalize \eqref{c>0}. 
 Using  results of section \eqref{General case} along with remarks  \eqref{Erlangrem},\eqref{Erlangrem2} we obtain

 \begin{theorem}   
    Suppose  A8 holds where  $f_{\tau_1}$ is given by  \eqref{hypoexp}. Let $ Q(s)\eq \prod_{  j}  (\lambda_j+s)$ and $ Q_0 = \prod_{  j}  \lambda_j $. 
 Then,       for $  \ 0\le x<b$, $N_b(x)$ solves    the integro-differential equation  

\beq     Q(-c\p_x )  N  (x)=   Q_0 \int_{0}^{x}  N(x-
 z)   h(z)  dz,     \label{Qeq}\ \end{equation} 
 with  BCs at the end-point $x=b$
  \beq  \partial_x^{j}N(b)= \delta_{j0},\ j=0 \dots n-1 \label{BC} \enq
 In particular, if $f_{\tau_1}(t)= \lambda ( \lambda t)^{n-1} e^{-\lambda t}/(n-1)!\sim \Gamma (n,\lambda)$ for some  $n\in\Bbb N,\lambda>0$ 
\beq    \Big(-c\p_x+\lambda\Big)^n N  (x)=  \lambda^n  \int_{0}^{x}  N(x-
 z)   h(z)  dz, \ 0\le x<b,  \label{Erlang} \end{equation}

\end{theorem}
  
   \begin{rem} \rm Note that -unlike \eqref{Fredac2}-  \eqref{Qeq}  has by itself  not a unique solution so  appropriate BCs  are required to pin down the EP.  Eqs. \eqref{Qeq} and \eqref{BC} define  a    final value problem  which  needs not be  well posed.      Section 5  elaborates on their derivation     under a fairly  general framework (see  \eqref{diffprob} and remark \eqref{Erlangrem}). 
\end{rem}   

   We next  {\it construct in explicit form} the solution to   \eqref{Qeq}.    To this end  we consider  {\it     the extension  from $(0,b)$ to   $ \Bbb R^+$ and     deprive it  of boundary conditions}.     Using  
  the known properties of  Laplace transformation $\mathcal L$:
 $$  \mathcal L(
   \Upsilon^{(n)})(s)=  s^n \hat \Upsilon (s)-  \sum_{j=0}^{n-1}  s^j    \Upsilon^{n-j-1}_0     \text{ and } \mathcal L(h\ast \Upsilon)(s)=\hat h(s) \hat \Upsilon(s)$$
 we obtain that 
any  continuous 
solution  $\Upsilon(x)$ of      \eqref{Erlang} on  $ \Bbb R^+$ must satisfy 
$$ \Big( Q (-cs )- Q(0)\hat h(s)  \Big)\hat \Upsilon(s)= \sum_{k=1}^n          (-c)^k  a_k   \sum_{j=0}^{k-1}  s^j    \Upsilon^{k-j-1}_0  $$ 
 $$ =\sum_{j=0}^{n-1}   s^j \sum_{k=j+1}^n           (-c)^k  a_k \Upsilon^{k-j-1}_0   \eq \sum_{j=0}^{n-1}   s^j\alpha_j(n)$$ 
 \beq \text{ where } Q(s)\eq  \sum_{j=0}^{n} a_j s^j, \  \alpha_j:= \sum_{k=j+1}^n           (-c)^k  a_k \Upsilon^{k-j-1}_0 \text{ and } \Upsilon^{(k) }_0\eq \p_x^k 
\Upsilon  (0)    \label{constant} \enq 

 Hence
\beq  \hat \Upsilon(s)=   \sum_{j=0}^{n-1}    \alpha_j s^j/ \Big( Q (-cs )- Q(0)\hat h(s)  \Big) \text{  }   \enq
  By inversion  we have that $ \Upsilon(x)$ can be written   in terms of $n$ arbitrary constants $\alpha_j $
   and a  {\it fundamental solution} $\pi(x)$ as:  
 
\beq  \Upsilon(x)=   \sum_{j=0}^{n-1}    \alpha_j \p_x^{j}\pi(x)  \quad \text{ where }  \    \pi (x)= \frac{1}{2\pi i}\int_{\beta-i\infty}^{\beta+i\infty} 
\frac{e^{sx}ds}{    Q(-cs)-Q(0)\hat h(s)  } \label{Erlanggennpi} \enq 

 
 Note that it  can be written in the suggestive way (compare with \eqref{Solfav2})

\beq     \pi (x)=  \frac{ 1}{2\pi i  }\int_{\beta-i\infty}^{\beta+i\infty} 
\frac{\bar f(-cs)}{ 1- \bar f(-cs)\hat h(s) }  e^{sx}ds\label{Erlangnpi2} \enq 

 \medskip 
  
  We have obtained a bundle of solutions $ \Upsilon(x)$    parametrized by   initial values  $   \p_x^k  \Upsilon  (0)$. The EP {\it should } follow     by imposing    the BCs \eqref{BC} for the values  $   N^{(k)}(b)$ at $x=b$. It is unclear that  the procedure works as this    problem needs not be   well-posed. We now prove that this is indeed  the case. 
Evaluation  of \eqref{Erlanggennpi} and its derivatives at $x=b$ implies  that the   constants $\alpha_j, j=0,\dots,n-1 $ must  satisfy  the linear system  $\mathbf{A}(b,b)\vec\alpha=  \vec 1$ where   $\vec 1=(1,0,\dots, 0)^\dagger$: 
 
\beq \mathbf{A}(b,b)\vec\alpha\eq \left(\begin{matrix}\pi& \pi^{(  1)} &\dots&\pi^{(  n-1)} \\
\pi^{(  1)}&  \pi^{(  2)}&\dots&\pi^{(  n)} \\\dots&&\dots\\  \pi^{(  n-1)}& \pi^{(  n)}&\dots&\pi^{(2n-2)}\end{matrix}\right) \left(\begin{matrix}\alpha_0 \\
 \alpha_1\\\dots \\ \alpha_{n-1} \end{matrix}\right)= \left(\begin{matrix} 1 \\
 0\\\dots \\ 0 \end{matrix}\right)\label{mat}\enq  
 Here  all matrix elements  $\mathbf{a}_{ij}=\pi^{( i+ j)}\eq   \p_x^{i+ j}\pi(x=b), i,j=0,\dots,n-1$
  are evaluated at $x=b$.  Besides   the Wronskian  of   the functions $\pi, \pi^{(  1)},\dots\pi^{(  n-1)}$ at $x=b$ is   \beq\Delta(b,b)= \det\mathbf{A}(b,b)\eq W\Big(\pi, \pi^{(  1)},\dots\pi^{(  n-1)}\Big) (x=b) \label{Wronski}\enq
 
   Let $\delta_j(b)$  be  the determinant of the matrix obtained substituting the $j-$th column of the matrix $\mathbf{A}(b,b) $   by the column vector $\vec 1$.    Cramer's rule  gives

  \beq \alpha_j=\frac{\delta_j(b)}{ \det\mathbf{A}(b,b)} \text{ and }
\  N_b(x) =   \underset{j=0}{\overset{n-1}\sum} \frac{\delta_j(b)}{ \det\mathbf{A}(b,b)}
 \p_x^{j}\pi(x)     \label{alphaN}  \enq        We introduce the $n\times n$ matrix  $\mathbf{A}(x,b)$ and  $\Delta(x,b)\eq \det \mathbf{A}(x,b)$ via

\beq \mathbf{A}(x,b)= \left(\begin{matrix}\pi(x)& \pi^{(  1)}(x) &\dots&\pi^{(  n-1)}(x) \\\pi^{(  1)}(b)&  \pi^{(  2)}(b)&\dots&\pi^{(  n)}(b) \\\dots&&\dots\\  \pi^{(  n-1)}(b)& \pi^{(  n)}(b)&\dots&\pi^{( 2n-2)}(b)\end{matrix}\right) \label{Amat}\enq

Note   that here  all  matrix elements   are evaluated at $x=b$   except for   those at the  first row. 
Let $\mathbf{A}_{0j}$ be  the $(0,j)$  minor of $\mathbf{A}(x,b)$, the determinant of the  matrix  that arises deleting the  $0-$  row and $j-$ th column. Then$ \mathbf{A}_{0j} =(-1)^{j} \delta_j$.  By row expansion  
   \beq   \det\mathbf{A}(x,b)=  \underset{j=0}{\overset{n-1}\sum}  (-1)^{j}
 \pi^{(j)}(x)   \mathbf{A}_{0j} (b) = \underset{j=0}{\overset{n-1}\sum}  
\delta_j (b)\pi^{(j)}(x)  \enq 

  Self-consistency of this procedure requires that $\p^{ j}\Delta(0,b)=\delta_{jn}  \Delta(b,b)$.      We skip   the proof which follows using $ \pi^{(j)}(0)=\delta_{j,n-1}$ . The following result summarizes the above.

 \begin{theorem}\label{5} Suppose     A8 holds and $f_{\tau_1}$ is given by  \eqref{hypoexp}(in particular, $\tau_1\sim \Gamma (n,\lambda)$ for some  $n\in\Bbb N,\lambda>0$). Then  

\begin{enumerate}\item The  integro-differential equation  \eqref{Qeq}  has  general solution  $ \Upsilon(x)=   \sum_{j=0}^{n-1}    \alpha_j \p_x^{j}\pi(x)  $  where   $\alpha_j, j=1,\dots n$  are arbitrary constants and  $\pi(x)$ is given by \eqref{Erlanggennpi} or  \eqref{Erlangnpi2}.\item  The  escape probability is given   in terms of   Wronskian determinants \eqref{Wronski} as 

 \beq  N_b(x)  \eq \Bbb P^x\Big(\tau^b<\tau^0\Big) = \frac{  \det \mathbf{A}(x,b)}    
{  \det \mathbf{A}(b,b)}   \quad    \label{escape4}   \enq

 \end{enumerate}
\begin{rem} The above result could be used to obtain survival probabilities by letting  $S(x)    =\underset{b\to\infty }\lim  \det \mathbf{A}(x,b) /  \det \mathbf{A}(b,b)$. This will be the subject of a future work. \end{rem}

\end{theorem}

   \section{  Risk model   under  rational   arrival  times} \label{General case}
Denote by   $\mathcal M$  the class of densities $f$ having rational    Laplace transform (LT)   $\hat f$:  \beq f\in\mathcal M\Leftrightarrow \hat f(s)=\frac{R(s)}{Q(s)}, s\ge 0 \label{Lapl} \enq    where  $Q,R$  are co-prime polynomials of orders $m\eq  \deg(R) <  \deg(Q)=n $:      \beq  \text{    }  Q(s)\eq \underset{j=0}{\overset{ n}\sum}a_j  s^j, \ R(s)\eq \underset{j=0}{\overset{ n-1}\sum}b_j  s^j \label{Lapl2}  \enq
     The characterization of such class   is not straightforward: a criteria in in terms of  complete monotonicity and unimodality was given  by Feller \cite{Feller3} and Bernstein; this approach is pursued in  \cite{Sumita}).   Obviously $Q(0)=R(0)\eq a_0\ne 0$ and   roots of $Q$ {\it must be located in the   negative real axis}. Besides, with no loss of generality, $a_n =1$. 

  We now establish several  results   that   relate  $\mathcal M$   with  solutions of   certain ordinary differential equations (ODEs).   The proof is deferred to  appendix B. 
\begin{Lemma}\label{Lemma1}    
 A density $f \in   \mathcal M$  
    iff  it  is of class
$C^n$ on $(0,\infty)$,  $n\in\Bbb N$ and   solves the ODE
 \begin{equation} Q(\p_t)  f \eq   \Big(\underset{j=0}{\overset{ n}\sum}a_j  \f {\p^j \ }{\p t^j}\Big)
f   =0 \ \text{ and  }\label{Geneq}\enq    where the   initial  data $ f_0^{(k)}\eq \p_x^kf(0), k=0,\dots n-1$  solve the linear system  \beq  \sum_{k=0}^{n-j-1} a_{j+k+1}    f^{(k)}_0=b_j, j=0,1,\dots n-1 \label{linear}\enq
 
\end{Lemma}

\begin{Corollary}\rm \rm    It follows from  \eqref{Geneq}    that 
 $  a_0\bar F(t)=\sum_{j=1}^n a_j \p_t^{ j-1}f(t)\label{Geneq2}$  
 \end{Corollary}
  
   ``Vectors''    $\vec f^{}_0\eq ( f^{0}_0,\dots,  f^{ (n-1)}_0)$ and $\vec b \eq( b_0,\dots, b_{ n-1})$  and    integer $m:=Deg   R$  have   a direct bearing on the degree of complexity of  eq.  \eqref{diffprob}  below which governs EPs.  Here we analyze their structure.  

\begin{Lemma}  \label{Er}  Let $ j_*\eq \min \{j:    f ^{(j)}_0\ne 0\}  $.      If $  f ^{(j)}_0=0$  for some $j, 0\le j\le n-2$      then  
\beq   f^{(k)}_0=0 \text{  for all } k\le j \text{ and } b_k=0  \text{  for all } k \ge n-j-1\enq   
   \beq  \text{ and }  j_* =  n-m-1    \label{m}\enq 
Thus $\vec f^{}_0 $ and $\vec b  $ have  at least one non-vanishing component $ f^{(n-1)}_0 \ne 0$ and {\it must have } the  structure (here    $\star$ denote a     non-null component)$$ \vec b^{} =( 
  \star,\overset{m+1}\dots \star,0,\dots 0), \  \vec f^{}_0 =(0,\overset{}\dots 0, \star,\overset{m+1}\dots \star)$$ 
    \end{Lemma}
   
\begin{rem}\label{Erlangrem}\rm  The above allows 
a partial classification of densities $f\in\mathcal M$   in terms of  the integers $n,m$ and $\nu_n,\nu_m$:  the {\it  number of different  roots of $Q$} and $R$ where 
\beq  1\le \nu_n\le n<\infty, \   1\le \nu_m\le m <n  \text{ or }\nu_m=0 \text{ if }m=0\label{restriction}\enq   Thus, for given $n$  a total of $n\Big(n^2-n+2\Big)/2$ sub-cases appear.  Some light is shed looking at the extreme cases: 
\begin{enumerate}\item $m=0$. This is the     {\it  hypoexponential } distribution previously studied. Besides it is   Erlang when  $  \nu_n=1$.   Here $\vec b^{} =(a_0, 0,
   \overset{n-1}\dots  0)$ and  $\vec f^{}_0 =(0,\overset{n-1}\dots 0,   a_0)$.   
\item  $m=\nu_m=n-1$ and $\nu_n=n$. Feller (\cite{Feller3}, pp. 439) proves  that   this corresponds to  {\it convex    mixture of exponentials} under the additional condition  
  $$ 0<\lambda_1<\beta_1<\lambda_2\dots  <\lambda_{n-1}<\beta_{n-1}<\lambda_n   \text{ where }$$ 
\beq   Q(s)=(s+\lambda_1)\dots (s+\lambda_n),  R(s)=(s+\beta_1)\dots (s+\beta_{n-1})  \enq
  
  \end{enumerate}
\end{rem}
 

 
  \subsection{{\bf Escape probabilities under      arrivals with  rational LT}}\label{rationalLT}
  We now  study EPs for the risk model    when A8 holds and       $f\in\mathcal M$. Such general case  is far more involved but can still be solved analytically by appropriately  transforming     (\ref{Fredac2})  into something  amenable to  Laplace transformation.    
  
\begin{theorem}\label{rational}

Suppose   that   $c>0$, and assumptions 8 and \eqref{Lapl} hold. 
 Let   $\Bbb Q$ and  $\Bbb R$ be the differential operators $ \Bbb Q\eq Q(-c \p_x), \Bbb R\eq R(-c \p_x)  $ and 
\beq q(x)\eq \int_{0}^{x} dzN( x-
 z)  h(z) dz \eq h   \ast  N \label{q}\enq  

  \begin{enumerate} 
\item The  solution $  N(x)   $ 
 of the integral eq. (\ref{Fredac2})  is       of class $C^n(0,b)$  and   must
 also   solve     on $0\le x<b$ the integro-differential   equation $ \Bbb Q N -\Bbb R ( q)=0$, or 
 
         \medskip

  \begin{equation}   \underset{j=0}{\overset{ n}\sum}(-c)^ja_j   \p_x^j N  (x)= \sum_{j=0}^{n-1}     (-c)^jb_j\p_x^jq(x)    \  \label{diffprob} \end{equation}
  \medskip

  and the $n-$BCs of terminal type (we denote $\xi_k\eq  (-1/c)^{k+1}  f^{(k)}_0,    f_0^{(k)}\eq \p_x^kf(0)$)
$$ N(b)=1, N'(b)=-\xi_0(1-q(b)),\dots $$ \beq \partial_x^{j}N(b)=   \sum_{k=0}^{j-1} 
 \xi_k q^{(j-k-1)}(b)-\xi_{j-1}, \ j=1 \dots n-1 \label{BCG} \enq 

\item Let    $\alpha_j, \ j=0 \dots n-1 $  be     free constants and $L$ be  the function   
\beq L(s)=  Q(-cs)-   R(-cs) \hat h(s)  =R(-cs)  \Big(1- \hat f(-cs) \hat h(s)\Big)/ \hat f(-cs) \label{G}\enq  Then a general solution to \eqref{diffprob}  is
\beq  \Upsilon(x)=   \sum_{j=0}^{n-1}     \alpha_j \p_x^{j}\pi(x)   \text{ where} \   \pi (x)= \frac{1}{2\pi i}\int_{\beta-i\infty}^{\beta+i\infty} 
      \frac{e^{sx} }{    L(s) } ds  \label{Genpi} \quad \text{} \enq 
   \end{enumerate}
 \end{theorem}

 {\it Proof}.     Operating with $\p^j_x$ on   \eqref{Fredac2}   we find,
  for $  j=1,\dots n$
\beq  \p_xN(x)=  (1/c) \Big( f (t_b)- f_0q(x)- (1/c)  \int_x^b q(z)    f^{(j)}((z-x)/c)dz  \Big ), \dots   \enq

 \beq  \p^j_xN(x)=(-1/c)^{j+1}\int_x^b q(z)   f^{(j)} ((z-x)/c)dz  +$$ $$ 
   +\xi_0 q^{ (j-1) }(x) +\xi_1 q^{ (j-2) }(x) +\dots \xi_{j-1}  q^{ }(x)  - (-1/c)^{ j }  f^{ (j-1) }(t_b)   \label{68} \enq

 Letting $x \to b^-$ the boundary conditions follow. 

With appropriate arrangement of the resulting terms we find after  some lengthy   calculations  
$$  \Big(  \underset{j=0}{\overset{ n}\sum}(-c)^ja_j   \p^j_x\Big)  N  (x)=E_1+E_2+E_3$$
$$ \text{where }\ E_1:=(1/c) \int_x^b dz q(z)    \sum_{j=0}^n   a_j   f^{ (j)}((x-z)/c)   \text{ and }$$
$$E_2:=   a_0\bar F((b-x)/c)-\sum_{j=1}^n a_j  f^{(j-1)}((b-x)/c) $$
 Eq. \eqref{Geneq} implies that the first two terms vanish: Concretely, $   Q(\p_t)  f =0$  yields  $E_1=E_2=0$.
   Upon simplification and using  \eqref{linear}  the third term is, 
$$E_3:= \sum_{j=0}^n (-c)^j a_j \sum_{k=0}^{j-1}(-1)^{k+1}  f^{(k)}_0  q^{( j-k-1)}(x)=$$$$ \sum_{k=0}^{n-1}\sum_{m=0}^{n-k-1} (-c)^m a_{m+k+1}    f^{(k)}_0 \p_x^m q(x)  = \sum_{m=0}^{n-1}   (-c)^m b_m \p_x^mq(x)   \qquad       $$

 \medskip

Hence $\Bbb QN(x)=\sum_{m=0}^{n-1}   (-c)^m b_m \p_x^mq(x) $ and     \eqref{diffprob} follows. 

 We solve   {\it  an auxiliary version of   \eqref{diffprob}   extended to   $0\le x<\infty$   deprived of boundary conditions}.   Laplace transformation      yields that any  
solution  $\Upsilon(x)$   must satisfy

$$Q(-cs)\hat \Upsilon(s) - \sum_{k=0}^{n-1}   u_ks^k= R(-cs)\hat q(s) - \sum_{k=0}^{n-2}    \eta_ks^k$$
where we introduce 
$$   \eta_k= \sum_{j=k+1}^{n-1}(-c)^j b_j  q_0^{(j-k-1)},  \  \eta_{n-1}=0    \text{ and }  u_k=\sum_{j=k+1}^{n} (-c)^j a_j  \Upsilon_0^{(j-k-1)} ,k=0,\dots, n-1 $$ 
The initial values  $  \Upsilon_0^{(j-k-1)}\eq \p_x^{(j-k-1)} \Upsilon(0) $ and $ q_0^{(j-k-1)}\eq \p_x^{(j-k-1)}q(0)$ are undefined so far. 
 It follows that 
$$ \hat\Upsilon(s)=\sum_{k=0}^{n-1}   (u_k-\eta_k) \frac{s^k}{L(s)}   $$


   By inversion we find the general solution    \eqref{Genpi} $ \quad \square$

\medskip

   In the  general case $ m\ge 0$ obtention of  the  EP   is far more involved than that of section 4.  We now  work the details.

\begin{theorem} \label{Gensol}   
Suppose   that assumption 8 and   \eqref{Lapl} hold.  Define 
  \beq  m_0(x)  = \int_{0}^{x}  \pi^{ } (z)    h(x-z)dz, \quad  m_j(x)  = \int_{0}^{x}  \pi^{(  j)} (z)    h(x-z)dz\label{m}\enq  
   Recall $\xi_k\eq  (-1/c)^{k+1}  f^{(k)}_0,k\ge 0$ and $\xi_{-1}:=-1$. Let        $\mathbf{A}=( \mathbf{a}_{ji}),  i,\ j=0 \dots n-1$ and  $   (\mathbf{\Theta} (x,b)   $  be  the $n\times n$  (respectively, $(n+1)\times (n+1)$) matrices  with entries   
\beq      \mathbf{a}_{0i} =
 \pi^{(  i )} (b), \quad   \mathbf{a}_{ji} =
 \pi^{(  i+j)} (b) -  \sum_{k=0}^{j-1}  \xi_k m_i^{(j-k-1)}(b),  i,\ j=0 \dots n-1, \label{a}  \enq   
\beq  {\bm \Theta} (x,b)= \left(\begin{matrix}0& \pi(x)  & \pi'(x)&\dots &\pi^{(n-1)}(x)
 \\ \xi_{-1}&   \mathbf{a}_{00}& \mathbf{a}_{0,1}  &\dots  & \mathbf{a}_{0,n-1} \\
  \xi_0 &        \mathbf{a}_{10} &        \mathbf{a}_{11}  &        \dots& \mathbf{a}_{1,n-1} \\  \vdots& &&& \vdots \\   \xi_{n-2}  &   \mathbf{a}_{n-1,0} &   \mathbf{a}_{n-1,1}  &        \dots& \mathbf{a}_{n-1,n-1}  \end{matrix}\right)_{(n+1)\times (n+1)}\label{Theta} \enq  

 Then the  EP is  
 \beq  N_b(x)    = \frac{   \det {\bm \Theta(x,b)}}    
{    \det {\bm \Theta}(b,b) }     \label{escape6}   \enq  \end{theorem} 

\begin{rem}\label{Erlangrem2}\rm    When $m=0$ the equation for the {\it  EP  \eqref{diffprob}    simplifies to \eqref{Qeq}}.   

Additionally  Lemma \eqref{Er} and \eqref{m} give $ \xi_j=0, 0\le  j  <n-1$; hence    $\mathbf{a}_{ji} =
 \pi^{(  i+j)} (b) $ and       all   entries but one   of the first column of  \eqref{Theta}  vanish.        Besides 
 $$a_0=Q(0)=R(0) =  \lambda_1\dots \lambda_n=  b_0=\sum_{k=0}^{n-1} a_{k+1}    f^{(k)}_0=  f^{(n-1)}_0 $$     
 giving    $ L(s)=  Q(-cs)-   Q(0) \hat h(s) $ and $    {\bm \Theta} (x,b)  $ coincides with \eqref{Amat}. 

\end{rem}
{\it Proof}.    Require  \eqref{diffprob}    to satisfy \eqref{BCG}.  Note  
$$  N(x)\eq  \sum_{i=0}^{n-1}     \alpha_i \pi^{(i)} (x) \text{ and } q(x)= \sum_{i=0}^{n-1}     \alpha_i  \int_{0}^{x}  \pi^{(i)} (z)    h(b-z)dz= \sum_{i=0}^{n-1}     \alpha_i   m_i(x),   $$ 
$$ N(b)\eq \sum_{i=0}^{n-1}     \alpha_i \pi^{(i)}(b)=1, \quad   N'(b) -\xi_0q(b) \eq  \sum_{i=0}^{n-1}  \Big(\pi^{(i+1)}(b)-\xi_0  m_i(b)\Big)\alpha_i =-\xi_0$$

     More generally,  it follows from \eqref{BCG} that (we denote $ m_i^{(j-k-1)}\eq \p_x^{( j-k-1)}m_i $)
  $$   \sum_{k=0}^{j-1} 
 \xi_k q^{(j-k-1)}(b)     =  \sum_{k=0}^{j-1} \xi_k  \sum_{i=0}^{n-1}    \alpha_i    m_i^{(j-k-1)}(b)  = \sum_{i=0}^{n-1}      \mathbf{c}_{ji} \alpha_i,  \ j=1 \dots n-1    $$ 
 $$\text{and } \ \partial_x^{j}N(b)\eq  \sum_{i=0}^{n-1}    \alpha_i    \pi^{(  i+j)} (b)= \sum_{i=0}^{n-1}      \mathbf{c}_{ji} \alpha_i   -\xi_{j-1},  \text{ or } $$
where at this stage we introduce the   the $n\times n$  matrice    $\mathbf{a}=( \mathbf{a}_{ji}), \mathbf{c}=(\mathbf{c}_{ji})$  with entries   
\beq \mathbf{a}_{ji} =
 \pi^{(  i+j)} (b) - \mathbf{c}_{ji} \text{ where }   \mathbf{c}_{0i}=0,   \mathbf{c}_{ji}=  \sum_{k=0}^{j-1}  \xi_k m_i^{(j-k-1)}(b), j\ge 1 \enq 
  Defining   $\xi_{-1}\eq f_0^{(-1)}\eq -1$   we  have $  N(x)\eq  \sum_{i=0}^{n-1}     \alpha_i \pi^{(i)}  $ where   $\vec\alpha  =(\alpha_0,\dots,\alpha_{n-1})$  solves  
\beq  \sum_{i=0}^{n-1}     \mathbf{a}_{ji} \alpha_i   =   -\xi_{j-1}\label{lineareq}, j=0,1\dots n-1  \enq 
  
Actually, a good deal more can be said about  the solution: By linearity  one has  
\beq \vec \alpha= - \sum_{k=0}^{n-1}   \xi_{k-1}   \vec \alpha^{(k)} \label{alpha^}\enq 
 where $\vec \alpha^{(k)}\eq   ( \alpha_0^{(k)}, \dots, \alpha_{n-1}^{(k)}), k=0,\dots n-1 $ solves the system 
\beq  \sum_{i=0}^{n-1}     \mathbf{a}_{ji}  \alpha^{(k)}_i   = \delta_{kj} \enq 
Cramer's rule  yields that \beq  \alpha^{(k)}_j=  \det    {\bm \theta}_j^k  /  \det\mathbf{A}\text{ and }  \alpha_j= - (\sum_{k=0}^{n-1}   \xi_{k-1} \det    {\bm \theta}_j^k )/  \det\mathbf{A} \label{alphaCra} \enq  where $\mathbf{A}= (\mathbf{a})_{ij} $ and $ {\bm \theta}_j^k  $ is the $n\times n$   matrix  obtained substituting the $j-$th column of  $\mathbf{A}  $   by the column vector $ (\vec e_{k})_m=\delta_{km} $: 
$$ ({   \bm \theta}_j^k)_{mn}=\mathbf{a}_{ mn}(1-\delta_{nj})+\delta_{mk}   \delta_{nj}, \ 0\le m,n\le n-1  \text{ or }
$$  
 
 \beq     {\bm \theta}_j^k   = \left(\begin{matrix}
  \mathbf{a}_{00}&     \dots0&\dots&  \mathbf{a}_{0,n-1} \\ \dots\\ \mathbf{a}_{k0}&     \dots 1&\dots&  \mathbf{a}_{k,n-1}   \\& &\dots\\   \mathbf{a}_{n-1,0} &\dots     0&\dots  &\mathbf{a}_{n-1,n-1}  \end{matrix}\right)  \enq

 Call   $ \mathbf{A}_k$     the matrix that results when  the $k-^{th}$ row of $ \mathbf{A} $ is substituted by the vector  $  ( \pi(x),  \pi^{(  1)}(x), \dots\pi^{(  n-1)}(x)) $, namely  

  $$  \mathbf{A}_k(x,b) = \left(\begin{matrix}
  \mathbf{a}_{00}&     \dots&   \mathbf{a}_{0,n-1} \\ \dots\\ \pi(x)&    \dots&  \pi^{(n-1)}(x)  \\& &\dots\\   \mathbf{a}_{n-1,0}       &\dots  &\mathbf{a}_{n-1,n-1}  \end{matrix}\right) $$

    The Laplace  co-factor expansion of this determinant yields  that 
 \beq    \det {\bm  A}_k(x,b) = \sum_{j=0}^{n-1}      \pi^{(  j)} (x) \det    {\bm \theta}_j^k  \label{A_k}\enq 
 A similar co-factor expansion  of matrix \eqref{Theta} gives

 \beq  \det {\bm \Theta}(x,b) = -\sum_{k=0}^{n-1}   \xi_{k-1}  \det \mathbf{A}_{k}(x,b)\label{ThetaA_k} \enq  
 
 Hence \eqref{alpha^}-\eqref{ThetaA_k}    yield  
$$  N_b(x)= \sum_{j=0}^{n-1}     \alpha_j  \pi^{(j)}(x)  = - \sum_{k=0}^{n-1}    \xi_{k-1}   \sum_{j=0}^{n-1}      \pi^{(  j)} (x) \vec \alpha^{(k)}_j $$ 
$$ = (-1/\det\mathbf{A}) \sum_{k=0}^{n-1}  \xi_{k-1}     \sum_{j=0}^{n-1}      \pi^{(  j)} (x) \det     {\bm \theta}_j^k   =     (-1/\det\mathbf{A}) \sum_{k=0}^{n-1}  \xi_{k-1}         \det \mathbf{A}_k =\frac{\det {\bm \Theta}(x,b)  }{ \mathbf{A} (b) }$$ 
 

  We next prove that  $\det{\bm \Theta} (b,b) =\det \mathbf{A} (b) $. Indeed,  $\det{\bm \Theta} (b,b) =$

 $$\begin{vmatrix} 0& \pi(b)  & \pi'(b)&\dots &\pi^{(n-1)}(b)
 \\ \xi_{-1} & \pi(b)  & \pi'(b)&\dots &\pi^{(n-1)}(b) \\
  \xi_0 &        \mathbf{a}_{10} &        \mathbf{a}_{11}  &        \dots& \mathbf{a}_{1,n-1} \\ \xi_1 &        \mathbf{a}_{20} &        \mathbf{a}_{21}  &        \dots& \mathbf{a}_{2,n-1} \\    \vdots& &&& \vdots \\   \xi_{n-2}  &   \mathbf{a}_{n-1,0} &   \mathbf{a}_{n-1,1}  &        \dots& \mathbf{a}_{n-1,n-1}   \end{vmatrix}=  \begin{vmatrix}    \pi(b)& \pi'(b)  &\dots  & \pi^{n-1}(b) \\
          \mathbf{a}_{1,0} &       \mathbf{a}_{1,1}  &        \dots& \mathbf{a}_{1,n-1} \\         \mathbf{a}_{20} &        \mathbf{a}_{21}  &        \dots& \mathbf{a}_{2,n-1} \\    \vdots& && \vdots \\     \mathbf{a}_{n-1,0} &   \mathbf{a}_{n-1,1}  &        \dots& \mathbf{a}_{n-1,n-1}  \end{vmatrix}= \det \mathbf{A} (b) $$

\subsection{ Escape probabilities when $n\eq Deg Q=2$} We use the former results to give {\it explicit expressions}  of EP when      Deg $Q\eq n=2$.  Let $0<p<1, q\eq 1-p$.     We have the  cases  (see remark \eqref{Erlangrem}): 
\begin{enumerate}
\item  $ m=0$,  $  \nu_n=1$ ({\it Erlang distribution}): 
   $\hat f(s)= \lambda^2/(\lambda+s)^2; \   f (t)=   \lambda^2 t   e^{-\lambda t} $. 
 \item   $ m=0$,  $  \nu_n=2$  ({\it hypoexponential distribution)}:
\beq \hat f(s)= \frac{ \lambda_1\lambda_2}{(\lambda_1+s) (\lambda_2+s) };0<\lambda_1<\lambda_2\Rightarrow \   f (t)=     \frac{ \lambda_1\lambda_2}{ \lambda_2-\lambda_1}\Big(   e^{-\lambda_1 t} - e^{-\lambda_2 t}\Big) \label{hypo} \enq  
 
\item $m= \nu_n=1 $  (Mixture of exponential and Erlang):     
\beq  \hat f(s)=   \lambda (ps +  \lambda)/(\lambda+s)^2 
 \Rightarrow \   f (t)=         \lambda  (p +  \lambda q t)e^{-\lambda t}     \enq   
\item   $  m=1, \nu_n=2 $ ({\it   Convex Mixture or hyperexponential}): 
\beq  \hat f(s)= \frac{(    p  \lambda_1+  q  \lambda_2) s + \lambda_1 \lambda_2 }{(s+  \lambda_1)(s+  \lambda_2)} 
 \Rightarrow \   f (t)=        p  \lambda_1   e^{-\lambda_1 t} + q  \lambda_2   e^{-\lambda_2 t} \label{mixture}  \enq  

 \end{enumerate}

\begin{exmp}{\it Escape probability under  hypoexponential and  $\Gamma(2,\lambda)$  distributions} \rm: 

  Suppose  $f_{\tau_1}$ is given by \eqref{hypo} with $ \lambda_1\le \lambda_2$.  Since  $m=0$ theorem \eqref{5} , \eqref{escape4}, \eqref{Erlanggennpi}  give 
 \beq  N_b(x)=\frac{\pi(x)\pi^{''}(b)- \pi^{\prime}(x)\pi^{\prime}(b)}{\pi(b)\pi^{''}(b)- (\pi^{\prime}(b))^2 } \text{ where}  \enq 
 \beq     \pi (x)= \frac{1}{2\pi i}\int_{\beta-i\infty}^{\beta+i\infty} 
\frac{e^{sx}ds}{     (\lambda_1-cs) (\lambda_2-cs)-\lambda_1\lambda_2 \hat h(s)  }  \enq

\end{exmp}

\begin{exmp}\rm Hyper-exponential:   $f (t)=        p  \lambda_1   e^{-\lambda_1 t} + q  \lambda_2   e^{-\lambda_2 t}$. 

  Hence  we suppose that $\hat f(s)$  is given by \eqref{mixture}  
where and $0 < p < 1, q=1-p$. Since $m=1$ the situation is considerably more complex, and the full formalism of theorem \eqref{Gensol} is required;  hence    we assume    $ -J_1\sim \mathcal E(\gamma)$. Since $   R(s)=  (p  \lambda_1+  q  \lambda_2) s + \lambda_1 \lambda_2,   Q(s)=(s+  \lambda_1)(s+  \lambda_2) $ the IE for $N(x)$ reads    (see \eqref{q}, \eqref{diffprob})

\begin{equation}     (\lambda_1-c\p_x) (\lambda_2-c\p_x)  N  (x)=      \Big( \lambda_1 \lambda_2-c (\lambda_1 p+\lambda_2 q)\p_x \Big)(h\ast N)  \end{equation} 
 Define $ \tilde \lambda_1\eq  \lambda_1/c , \tilde \lambda_2\eq  \lambda_2/c $.   The function $\pi(x)$ is retrieved from  \eqref{Genpi} and  \eqref{G} where   
 $$L(s)= \frac{c^{2}s}{\gamma+s}\Big(s^{2}+  A s+ B\Big), \ A\eq  \gamma - \tilde\lambda_1-\tilde\lambda_2 , \quad B\eq \tilde\lambda_1\tilde\lambda_2- \gamma(q\tilde\lambda_1+p \tilde\lambda_2 ) $$
Thus $L(s)$ has poles on  $\mathcal R=\{0,s_-,s_+\}$ where we define  $ s_{\pm }=  \Big(-A\pm \sqrt{A^2-4 B}\Big)/2$. 
  It follows that   (see \eqref{m})  
 

 \begin{equation}
\pi(x)=  \frac{1 }{ c^2  s_+s_- }    \Big(\gamma+\frac{ 1 }{s_+-s_-} \Big( e^{s_+x} s_-(\gamma+s_+)- e^{s_-x} s_+(\gamma+s_-)\Big) \Big),   \enq 
 $$   m_0(x)  =    \frac{1 }{ c^2  s_+s_- } \Big(       1+\frac{ 1 }{s_+-s_-} \Big(   s_-e^{s_+x} -  s_+e^{s_-x}  \Big)  \Big)$$ 
$$\text{ and }   c^2 (s_+-s_-) m_1(x)  =        s_-e^{s_+x} -  s_+e^{s_-x}  $$ 

      Besides $f_0\eq    \lambda_1 p  + q\lambda_2$ . Recalling that $\xi_0=-f_0/c$ and $ \mathbf{a}_{10} =
 \pi^{\prime)} (b) -      \xi_0 m_0 (b), $ and setting $\mu\eq  c^5 (s_+-s_-)s_+s_-$  we find the EP via \eqref{escape6} and  \eqref{Theta} where
  
 \beq  \mu \det  {\bm \Theta} (x,b)= \mu  \begin{vmatrix}0& \pi(x)  & \pi'(x) 
 \\ -1&   \pi(b)& \pi'(b)   \\
  \xi_0 &         \mathbf{a}_{10}  &        \mathbf{a}'_{10}    \end{vmatrix}= s_+(f_0-c(s_++\gamma))  e^{s_+b}  \Big( (s_-+\gamma) e^{s_- x}-\gamma\Big)        -  $$ 
 $$ - s_- (f_0-c(s_-+\gamma)) e^{s_-b} \Big( (s_++\gamma) e^{s_+ x}-\gamma\Big)   \label{Theta2} \enq 
 

\end{exmp}
 
\subsection{Ideas on the case $n=3$}   Suppose  Deg $Q\eq n=3$. Bearing in mind the restrictions \eqref{restriction}  there are up to 12 possible cases labeled by $(\nu_n,m,\nu_m)$  which we do not attempt to classify.   Consider however  the  following interesting cases   $f_{1,2}\in\mathcal M$
\begin{enumerate}
\item  Let  $0<p<1$ and $\lambda>0$ and 
 $ \hat f_1(s)=   \lambda \Big( \lambda^2+ 2sp \lambda  + s^2p \Big)/(\lambda+s)^3  
       $
    \item Take now $\lambda>1,\alpha=\Big(1\pm\sqrt{\lambda^3-1}\Big)/\lambda$ and 
\beq  \hat f_2(s)=   \Big(2 - 2\alpha(\lambda+s)  + \alpha^2(\lambda+s)^2 \Big)/(\lambda+s)^3\enq   
\end{enumerate}
  
   The corresponding  densities $f_{1,2} $ have equal  integers $m=\nu_m=2, \nu_n=1$     ($R$  having a  pair of complex conjugate roots); nevertheless  they  are markedly different. Actually, 
\beq  f_1 (t)=         \lambda  (p +  \lambda^2 q t^2/2)e^{-\lambda t}\in\mathcal M\text{ and } f_2 (t)=          (t-\alpha)^2e^{-\lambda t}\in\mathcal M\enq

\section{Two-sided problems}   In this section we address the   situation where jumps may take {\it both signs}.  
{\it It is remarkable that   \eqref{fav2}  is still  solvable when only one of the conditions   $q_1=0$ or $p_1=0$  is required   and the remaining parameters  $p_1,p_2, q_1, q_2$       }   are arbitrary. Recall that $\star$ signifies a  non-null component.    Table 1 summarizes these results.
\subsection{{\bf Support  $H \subset(-\infty,0)\cup[b-x,\infty)$ or $\vec p=(\star,\star,0,\star)$.} } Here we show that the  ideas of sections 3 and 4 carry over  to the case    support  $H \subset(-\infty,0)\cup[b-x,\infty)$. We  consider the case when jumps are  hypoexponential,  which helps to keep the algebra tidy.   Generalization to arrivals $f\in\mathcal M$ is messy but  straightforward.

 \begin{theorem}    Suppose assumptions  A1-A7 hold with $\tau_n $  given by  \eqref{hypoexp}  and   let $ Q(s)\eq \prod_{  j}  (\lambda_j+s)$.    Suppose  support  $H \subset(-\infty,0)\cup[b-x,\infty)$ where  $  p=1-q= \Bbb P( J_1\ge b-x )$  and   $ dH_-(x)=h(-x) dx$. Then 
 \begin{enumerate} \item $N_b(x)$     satisfies the BCs \eqref{BC} and solves        for $  \ 0\le x<b$

\beq     Q(-c\p_x )  N  (x)=  Q_0\Big(p+ q  \int_{0}^{x}  N(x-
 z)   h(z)  dz   \label{Qeq2}\ \end{equation} 
 
  \item Let $L(s):= Q(-cs)-   q  Q_0 \hat h(s)  $.  Define the fundamental solutions \beq     \pi (x)= \frac{1}{2\pi i}\int_{\beta-i\infty}^{\beta+i\infty} 
      \frac{e^{sx}ds }{   L(s)   }  \text{ and } \pi_{-1} (x)= \frac{1}{2\pi i}\int_{\beta-i\infty}^{\beta+i\infty} 
      \frac{e^{sx}ds }{ s  L(s)   }  \label{Genpi2} \enq  
Let   $\mathbf{A}(x,b)$ be  the $n\times n$ matrix \eqref{Amat} and  $\mathbf{B}(x,b)$ be the $(n+1)\times (n+1)$ matrix 
 
$$  \mathbf{B}(x,b)= \left(\begin{matrix}  \pi^{(-1)}(x)  & \pi(x)&\dots &\pi^{(n-1)}(x) 
 \\    \pi^{(-1)}_b  & \pi_b&\dots &\pi^{(n-1)}_b
 \\     \pi^{  }_b  & \pi^{(1)}_b&\dots &\pi^{(n)}_b  
 \\    \vdots& &&  \vdots \\         \pi^{(  n-2)}_b& \pi^{(  n-1)}_b &\dots&\pi^{(2n-2)} _b \end{matrix}\right) _{ n+1\times n+1}     $$
where $\pi_b\equiv \pi(b),\dots$.Then the   escape probability is 

   \beq  N_b(x)  \eq \Bbb P^x\Big(\tau^b<\tau^0\Big) = \frac{  \det \mathbf{A}(x,b)+ p Q_0\det \mathbf{B}(x,b)  }    
{  \det \mathbf{A}(b,b)}   \quad    \label{escape5}   \enq

  \end{enumerate} 

 \end{theorem}  
{\it Proof}.  We start noting that  $\pi_{-1} (x)=\int_0^x  \pi  (y)dy\eq  \pi^{(-1)}(x)$. Besides eq.    \eqref{fav2}  reads    
 \beq N(x)=    q\bar F (t_b) +  p  + (q/c) \int_{x}^b   dl    f((l-x)/c ) \int_{0}^{ l} N( l-y) h(y) dy \label{Fredac3} \enq 
  
Eq. \eqref{Qeq2} follows with similar reasonings to that used as for \eqref{Qeq} and \eqref{diffprob}. Note that using the results in  remark \eqref{Erlangrem}, \eqref{68}   is modified to 
 $$  (1/q) \p^j_xN(x)= (-1/c)^{j+1}\Big(\int_x^b q(z)   f^{(j)} ((z-x)/c)dz  + f^{ (j-1) }(t_b)- \delta_{jn}Q_0 h\ast N \Big),  1\le j\le n    $$
Laplace transformation   of \eqref{Qeq2}    yields  now

$$\Big( Q(-cs)  -  q Q_0 \hat h(s) \Big)\Upsilon(s)  =\frac{pQ_0}{s}+ \sum_{j=0}^{n-1}   s^j\alpha_j(n) $$  
for certain free constants $\alpha_j, j=0,\dots n-1$. A  general solution   is
\beq  \Upsilon(x)=  pQ_0 \pi_{-1}(x)+  \sum_{j=0}^{n-1}     \alpha_j \p_x^{j}\pi(x) \label{generalsol3} \quad \text{} \enq 
  To obtain $N(x)$ we require for $k=0,\dots n-1$ 
$$pQ_0 \pi^{(k)}_{-1}(b)+   \sum_{j=0}^{n-1}     \alpha_j  \pi^{(k+n)}(b)=\delta_{k0},$$ 
 Introducing  $\alpha_{-1}:=pQ_0$ the above can be written  as  $ \Upsilon(x)=     \sum_{j=-1}^{n-1}     \alpha_j \p_x^{j}\pi(x)  $ where 
  the $\alpha$'s satisfy  the system of $(n+1) $  equations  
  
\beq  \left(\begin{matrix}  1 & 0 &\dots & 0& 0  \\ 
\pi_{-1}& \pi& \pi^{(  1)} &\dots&\pi^{(  n-1)} \\
 \pi_{0}&      \pi^{(  1)}&  \pi^{(  2)}&\dots&\pi^{(  n)} \\\dots&&\dots\\  
 \pi^{(  n-2)}& \pi^{(  n-1)}& \pi^{(  n)}&\dots&\pi^{(2n-2)}\end{matrix}\right) \left(\begin{matrix}\alpha_{-1} \\ \alpha_0 \\
 \alpha_1\\\dots \\ \alpha_{n-1} \end{matrix}\right)= \left(\begin{matrix}  pQ_0   \\\
  1  \\0 \\\vdots \\ 0 \end{matrix}\right)\label{mat}\enq   
 
This has the form \eqref{lineareq} where $\xi_{-1}:=pQ_0, \xi_0:=1,  \xi_j=0, j\ge 1$. Hence, repeating mutatis-mutandis the arguments in theorem \eqref{Gensol}   we may obtain the solution  given by eqs. \eqref{Theta}, \eqref{escape6}, where ${\bm  \Theta} (x,b)$ is now   the  $(n+2)\times (n+2)$ matrix    with determinant

 $$ \det {\bm \Theta} (x,b)= \det  \left(\begin{matrix}0& \pi^{(-1)}(x)  & \pi(x)&\dots &\pi^{(n-1)}(x)\\ - pQ_0  &  1  &  0&\dots &0
 \\   -1 & \pi^{(-1)}_b  & \pi_b&\dots &\pi^{(n-1)}_b
 \\   0& \pi^{  }_b  & \pi^{(1)}_b&\dots &\pi^{(n)}_b  
 \\    \vdots& &&& \vdots \\    0 &    \pi^{(  n-2)}& \pi^{(  n-1)} &\dots&\pi^{(2n-2)} 
\end{matrix}\right)_{(n+2)\times (n+2)} =  $$
$$pQ_0\begin{vmatrix}  \pi^{(-1)}(x)  & \pi(x)&\dots &\pi^{(n-1)}(x) 
 \\    \pi^{(-1)}_b  & \pi_b&\dots &\pi^{(n-1)}_b
 \\     \pi^{  }_b  & \pi^{(1)}_b&\dots &\pi^{(n)}_b  
 \\    \vdots& &&  \vdots \\         \pi^{(  n-2)}& \pi^{(  n-1)} &\dots&\pi^{(2n-2)} \end{vmatrix}  + \begin{vmatrix}   \pi(x)&\dots &\pi^{(n-1)}(x)\\    
     \pi^{(1)}_b&\dots &\pi^{(n)}_b  
 \\    \vdots& &   \vdots \\           \pi^{(  n-1)} &\dots&\pi^{(2n-2)} \end{vmatrix}=       $$
$$ \det \mathbf{A}(x,b) +pQ_0 \det \mathbf{B}(x,b)$$
Clearly $ \det \mathbf{B}(b,b)=0 $ and  $\det {\bm \Theta} (b,b)= \det \mathbf{A}(b,b) $ and the result follows. 
\begin{Corollary}      Suppose  that $\tau_1\sim\mathcal E(\lambda) $ and  support  $H \subset(-\infty,0)\cup[b-x,\infty)$ where   $ dH_-(x)=h(-x) dx$ and  $  p=1-q=  \Bbb P( J_1\ge b-x )$.   Then \eqref{escape2} is generalized to 
  \beq \pi  (x)\eq \frac{1}{2\pi i}\int_{\beta-i\infty}^{\beta+i\infty} ds \frac{\lambda e^{sx} }{\lambda   q \hat h(s) - \lambda+cs }   \text{ and } \enq

\beq N(x)= \frac{\pi(x)}{\pi(b)} + p\lambda \Big(  \pi^{-1}(x) \pi(b)-\pi(x) \pi^{-1}(b)\Big)/\pi(b)\enq  
\end{Corollary}


\medskip

\subsection{{\bf  support  $H \subset(-\infty,-b)\cup[0,\infty)$}:  $\vec p=(\star, 0,\star,\star)$   } 
We next   show that when {\it    $q_1=0$  the solution to   \eqref{fav2} can be given in closed way  }   for  general distribution  of jumps and severities.   
 
 \begin{theorem}\label{gentrivial} Suppose assumptions  A1-A7 hold with $\tau_n\overset{i.i.d} \sim F$ and   $J_n\overset{i.i.d}\sim  H$ with densities  $F'=f$ and $H_+'\eq h_+$. Suppose  that $q_1=0$, i.e.    supp.  $H\subset (-\infty,-b]\cup [0,\infty)$.       Let  $\pi(x)$ be the solution of the integral equation defined on   $0\le x <\infty $

  \beq   \pi(x) =  \bar F (x/c)+  (p/c)  \int_0^{x}  \bar   H_+( z )f  ((x-z)/c)dz +   
   (p/c)  \int_0^{x}   f  ((x-z)/c)dz\int_0^zh_+(y) \pi(y-z) dy
 \label{Q2} \enq   
 \begin{enumerate}\item  $\pi(x)$ has Laplace transform
  \beq \hat {\pi} (s)= 
  \Big(1-  q\hat  f (cs) \Big(1- p   \hat h_+(s)     \hat  f (cs)\Big)^{-1}   \Big)/s\label{Laplsol} \enq  
 \item    Suppose that $ F(b/c)  p_1 <1$.  Then $ N_b(x)$   is given by   
\beq   N_{b} (x)=    \pi( b-x) = 
    1 -\int_{\beta-i\infty}^{\beta+i\infty} \frac{ q\hat  f (cs)e^{(b-x)s}}{2\pi i s\Big(1- p   \hat h_+(s)     \hat  f (cs)\Big)} ds\label{Solfav2}\enq  
\end{enumerate} \end{theorem} 
 
{\it Proof}.    If $q_1=0$   \eqref{fav2}  reads for $ \ 0\le x\le b$:  
 \beq N(x)= \bar F (t_b)+   p \int_0^{t_b}  dF(l)  \Big( \bar   H_+(b-x-cl ) +  \int_{0}^{b-x-cl}     N(x+cl+y) dH_+(y)\Big),  \label{fav}    \enq

  where $ 1-q=p\eq p_1+p_2=\Bbb P(J_1>0)$.  Setting $f_c(x):=f(x/c)$, 
    \eqref{Q2} implies
 $$   \pi(x)  =    \bar  F_c(x)+       (p/c)  \int_0^{x}dz f_c  (x-z)  \Big(   \bar   H_+( z ) \int_0^zh_+(y) \pi(y-z)  dy\Big) $$
 Recalling that $\mathcal L(  \bar F )(s)=(1-\hat f )/s$ and   $\mathcal L(   f_c )(s)=c \hat f  (cs)$ and noting that the above is a repeated convolution we find 
$$\hat \pi(s) =(1/s)  \Big(1-q\hat f -p\hat f\hat h\Big)(cs) + p\big( \hat f\hat h\hat\pi\big)(cs)  $$ 
   \eqref{Laplsol} follows.  
For (ii)  note that $\Pr(J_1>x)= p \bar   H_+(x)$.
 The key idea is to introduce a new function via $\tilde  N (x)\eq  N(b-x)$. Eq. \eqref{fav} is transformed to
$$    \tilde  N (x) =  \bar F_c (x)+  p  \int_0^{x/c} dF (l)  \Big( \bar   H_+(x-cl )+
     \int_{0}^{x-cl}    dH_+(y)       N(b-x+cl+y\Big)  $$  
$$=  \bar F_c (x)+  (p/c)  \int_0^{x}dz   f_c   (x-z)  \Big(   \bar   H_+( z ) +   \int_0^zh_+(y)  \tilde N(y-z)dy\Big)  $$ 
where we used $y\to z=x-cl$ and  $N(b-z+y)= \tilde N(y-z)$. 
 Thus   $  \tilde N (x) $  
satisfies  for  $0\le x\le b$   the same equation as  $ \pi(x)$  does on  $[0,\infty)$, namely \eqref{Q2}.  Since  the bound \eqref{bound} guarantees   the existence of a {\it unique} solution  both
functions must be the same:    $\tilde N_b(x)=    \pi(x)$ for  $0\le x\le b \quad \square$

\medskip 

We next consider several examples      of interest with  $q_1=0$ and  $p=\Bbb P(J_1>0)$. 

  \begin{exmp}{\it  Exponential jumps.} \rm    Suppose   $\tau_1\sim \mathcal E(\lambda)$,  supp.  $H_-\subset (-\infty,-b]  $ with $q= \Bbb P(J_1<0) $ and  $ H_+\sim \mathcal E(\gamma)$.   \eqref{Solfav2} gives 
    \beq   N_b(x)=
     \frac{ 1}{  \gamma (s_+ -s_-)}  \Big( e^{s_- (b-x)} s_+ (s_-+\gamma)  - 
e^{s_+ (b-x) }s_- (s_++\gamma) \Big) 
 \enq  
where 
$$2s_{\pm}=  -(\rho+\gamma)  \pm  
  \sqrt{  (\rho+\gamma)^2 - 4 q  \rho \gamma} . $$

\end{exmp}
     
 \begin{exmp}{\it fixed magnitude jumps.}\rm   We consider the case     $\tau_1\sim \mathcal E(\lambda)$,   $q_1=0$  while positive jumps
have a {\it  fixed magnitude } $y_1$, i.e.  $\hat h_+(s)=e^{-s y_1}$ and $p=\Bbb P(J_1>0)$. We see that

$$N_b(x)=1 - q \rho  \int_0^{b-x}dy C(y,\rho,p) \text{ }$$ 
where     $ C(y,\rho,p) $ was evaluated in example 4.  
  Note that 
  $$  \rho  \int_0^{x}dy   (   \rho  (y-k y_1))^k  e^{  - \rho(y-k y_1)} /k!=1- \sum_{j=0}^k
 (  \rho (x-k y_1 ))^k  e^{  - \rho(x-k y_1)} /  j!  $$ 
 Hence, if we define $ \xi_k\eq  b-x-k y_1, \rho\eq \lambda/c \text{ and }  k_1\eq  \lfloor (b-x)/y_1\rfloor$ we finally have

 \beq  N_b(x)= 1- q \sum_{k=0}^{  k_1 } p^k \Big( 1- \sum_{j=0}^k
 (   \rho \xi_k)^j  e^{  - \rho   \xi_k} /  j!  \Big)\label{solfix}   \enq   
  
Note that 
$\underset{k\to\infty}\lim  \sum_{j=0}^k
 (   \rho \xi_k)^j  e^{  - \rho   \xi_k} /  j!  =1$.   Besides when $y_1$ is so large as    $  b-x <y_1$ then $0=k_1=k=j$ and  \eqref{solfix}  reduces to  \eqref{double}. 
Finally as $p\to 1$  then $ N_b(x) \to 1$ as expected, 
\end{exmp} 
 
It is interesting to compare the EP correspnding to  $\tau_1\sim\mathcal E(\lambda)$ and the cases where    (i)   $ H \sim$    Laplace$(0,\gamma)$,  (ii)  $ H_+\sim \mathcal E(\gamma)$, $ \Bbb P(J_1<0)= \Bbb P(J_1\le -b)=1-p$, and (iii)    $ H_+\sim \delta (y-y_0)$ and $ \Bbb P(J_1<0)= \Bbb P(J_1\le -b)=1-p$. 


\subsubsection{ Severities with    rational  characteristic function}  In the spirit of section 5.1 we  denote by   $\mathcal H$  the class of densities $h$ having rational       characteristic function (CF), namely     
\beq  h\in\mathcal H\Leftrightarrow  \tilde h (\omega):=  \int_{-\infty}^\infty  h(x) e^{ix\omega} d\omega= \frac{R(i\omega)}{Q(i \omega)},  \label{Four}   \enq   
\beq  \text{  Here }   Q(s)\eq \underset{j=0}{\overset{ n}\sum}a_j  s^j, \ R(s)\eq \underset{j=0}{\overset{ n}\sum}b_j  s^j\enq
 are  co-prime  polynomials     with   $  \deg(R) <  \deg(Q)=n $. Besides $(a_j), (b_j)\in \Bbb R$,  $a_0=b_0$ and $a_n\ne 0, b_n=0$.    
 It turns out that  {\it for severities  $h\in\mathcal H$ and also   $N_t\sim\mathcal P(\lambda t)$    \eqref{Poissongen} can be reduced further to  an    ordinary differential equation} (ODE)  with boundary conditions at $x=0$.  
 Interesting examples of such class include 
 \begin{enumerate}
 \item   Suppose  positive (negative) jumps are exponentially distributed  with means  $\gamma_{\pm}^{-1}$ and let     $p:= \Bbb P( J_1>0), q=  1-p$. Then   \beq   \tilde h(\omega)=p\frac{\gamma_+}{\gamma_+ -i\omega } +q\frac{\gamma_-}{\gamma_- +i\omega }, \label{three-parameter}  \enq 
\beq  h(x)= p\gamma_+e^{-\gamma_+x}\theta (x)+ q\gamma_-e^{\gamma_-x}\theta (-x) \text{ and }   H(x)=   q    e^{\gamma_-x}\theta (-x) + 1-pe^{-\gamma_+x}\theta (x)  \label{three-parameterf}  \enq 
 Such {\it double-exponential jump  models}  find application in  mathematical finance.
  It corresponds to the polynomials
\beq   \quad R(s)=\gamma_-\gamma_++( p\gamma_+- q\gamma_-)s, Q(s)=(\gamma_+- s) (\gamma_-+ s),  \label{three-parameterR}  \enq 
 
\item   $J_1 \sim$    Laplace$(0,\gamma)$  is recovered when  $p=1/2$, 
$\gamma_+=\gamma_-:=\gamma >0$.  
\item   The {\it variance gamma } distribution (VGD) is a widely used  model      in  stochastic finance.  If $ \s, \vartheta\in\Bbb R, n\in\Bbb N $ and $a^2=\Big(2+\frac{\vartheta^2}{\s^2}\Big)\s^{-2}$ it is given by     
\beq  \tilde h(\omega)= \Big( (\s^2/2)\omega^2+i\vartheta \omega +1 \Big)^{-n} \text{ and } h(x)=     C e^{-\vartheta/\s^2 x -a|x|}P|(x|)   \label{variance}\enq  
  where $P$ is a certain polynomial with degree $n-1$ and $C$ a normalizing constant.  
 \end{enumerate} 
We   first establish 
the following Lemma,  which is proved in  Appendix C.

 \begin{Lemma}      Assume  $  h\in\mathcal H $ with $\tilde h$  given by  \eqref{Four}. Let   $ I_j\eq h^{(j)}(0^+)-h^{(j)}(0^-)$   be   the jump at the origin of $\p_x^j h, j \le  n-1$.   
 Then $h(x)$ solves    the ODE  with boundary conditions at $x=0$
  \beq Q(-\p_x) h_{} (x)=  0, \ x\in \Bbb R-\{0\}\label{ODE}\enq \beq         \text{ and }  \underset{j=k+1}{\overset{ n}\sum} (-1)^{j-k}a_j I_{j-k-1}= b_k, \ k=0,\dots n-1    \label{ODEBC}\enq 
  
 Reciprocally if  $h$ is a density and solves  \eqref{ODE}, \eqref{ODEBC}  then it has a CF  given by 
 \eqref{Four} $\square$ 

\end{Lemma}

We now show that  EP can be found solving a simple ODE. 
 
     \begin{proposition}\label{Four2} Suppose assumptions A1-A7 hold with $\tau_1\sim\mathcal E(\lambda)$ and    $h\in\mathcal H$  satisfiying  \eqref{Four}.     Let  $\Bbb L$ be the differential  operator $\Bbb L\eq  \Big(Q-R -\rho^{-1} \p_x \circ Q\Big)(\p_x)$. Then $N_b(x)$  solves the ODE  
\begin{equation}  \Bbb LN\equiv  \Big(\underset{j=0}{\overset{ n}\sum}  (a_j-b_j) \f {\p^j \ }{\p x^j} -\frac{  a_j }{\rho}  \f {\p^{j+1}   }{\p x^{j+1}} \Big)N(x)=0\label{solvable2}\end{equation}

         Further,  $n_j\eq N^{(j )}(b^-), j=0,\dots n-1    $   satisfies  the  linear  system of BC$: n_0 =  1$ and 
$$    \ n_0- \rho^{-1}  n_1 = \bar H(0^+)+      \int_{0}^{b}  h_-(z-b)  N(z)dz,  \text{ and for }  j=1,\dots n-1   $$  
      \beq   n_j - \rho^{-1}  n_{j+1}    +\underset{k=0}{\overset{ j-1}\sum}(-1)^{k}I_k    n_{j-k-1}= (-1)^{j-1}\Big(h^{(j-1)}_+(0) - \int_{0}^{b}  \p_z^j h_- (z-b) N(z)dz\Big)\label{BC3} \enq  
\end{proposition}
 {\it Proof}.     We write the jump distribution as 
 $ h(y)=     h_{+}(y)\theta(y)+  h_{-}(-y)\theta(-y) $
where  $\theta(x)=\mathbf{1}_{x\in(0,\infty)}$, the Heaviside  function. Since  $\tau_1\sim\mathcal E(\lambda)$ then Eq. \eqref{Poissongen} applies. To keep the algebra  tidy we introduce      $M(x):= N(x)- \rho^{-1}\p_x N(x)$ and \eqref{Poissongen} reads
 $$  M(x)=   \bar H(b-x)+    \Big(  \int_{x}^{b} h_+ (z-x)+ \int_{0}^{x}  h_-(z-x)\Big) N(z)dz   $$ 
By repeated differentiation we find for $j\ge 1$
$$\p_x^{(j)} M(x)=  (-1)^{j-1}h^{(j-1)}_+(b-x) + \underset{k=0}{\overset{ j-1}\sum}(-1)^{k+1}I_k    N^{(j-k-1)}(x)+\Big(  \int_{x}^{b}  + \int_{0}^{x}  \Big)(-1)^j\p_z^j h (z-x) N(z)dz$$ 
  This  yields the BC \eqref{BC3} sending $x\to b$.   

Next, operating with $ Q(\p_x)$  on the LHS of \eqref{Poissongen}   yields that  $M$ satisfies 
$$ \sum_{j=0}^n a_j  \p_x^{(j)} M(x)= a_0    \bar H(b-x) +a_1 h_+(b-x) - a_2 h'_+(b-x) +\dots   (-1)^{n-1}h^{n-1}_+(b-x) +$$$$  \sum_{j=0}^n a_j \underset{k=0}{\overset{ j-1}\sum}(-1)^{k+1}I_k    N^{(j-k-1)}(x)+\Big(  \int_{x}^{b}  + \int_{0}^{x}  \Big) N(z)\sum_{j=0}^n a_j  (-1)^j\p_z^j h (z-x) dz$$
   Recalling \eqref{ODE} we see that several terms cancel as
 $\sum_{j=0}^n a_j  (-1)^j\p_z^j h (z-x)\eq Q(-\p_z)h(z-x)=0$.   The above simplifies to  
$$  \sum_{j=0}^n a_j  \p_x^{(j)} M(x)=\sum_{j=0}^n a_j \underset{k=0}{\overset{ j-1}\sum}(-1)^{k+1}I_k    N^{(j-k-1)}(x) =$$$$ \sum_{m=0}^{n-1}  N^{m)}(x)  \underset{j=m+1}{\overset{ n}\sum}(-1)^{j-m}I_{j-m-1}a_j  =  \sum_{m=0}^{n-1}  N^{m)}(x) b_m =R(\p_x)N$$ 
where we used \eqref{ODEBC}. 
 Eq. \eqref{solvable2} follows since 
$$   Q(\p_x)M\eq    Q(\p_x)N-  \rho^{-1}\p_x Q(\p_x)N = R(\p_x)N\quad \square $$

   We next evaluate the EP for  several cases of interest when  $\tau_1\sim\mathcal E(\lambda)$. 
\begin{exmp}{Risk model recovered. }\rm   To warm up suppose again  $-J_1\sim \mathcal E(\gamma)$ and $J_1<0$. This entails    (see\eqref{three-parameterR})
 $R(s)=\gamma, Q(s)=\gamma+ s $.  From \eqref{solvable2}  the EP is found solving  $$ \Big(    \p_{xx}+ (\gamma-\rho)\p_x\Big)N(x)=0, \text{ and } N(b)=1, \quad   \rho-   \rho   \int_{-b}^{0}  N (b+y)  e^{-\lambda y}dy=N'(b)$$
One checks easily that  \eqref{Poissexp} is the only solution to this ODE and  BCs.    \end{exmp}
 
\begin{rem}\rm  Notice that the   EP for   Brownian motion with drift $v\eq   (\gamma-\rho)/2$ {\it satisfies also}  $\Bbb L N=0$ and that its infinitesimal generator is, up to a constant, $\Bbb L\eq \p_{xx}+   (\gamma-\rho)\p_x$.   \end{rem} 

 \begin{exmp}{Laplace distribution. }\rm  Suppose that   $J_1 \sim$    Laplace$(0,\gamma)$. It follows that  (see  eq. \eqref{three-parameterR}) $ R(s)=\gamma^2, Q(s)=\gamma^2-s^2 $  
 and  $N_b(x)$ must satisfy 
$$  \Big( \p_{xxx}-\rho  \p_{xx}- \gamma^2\p_x\Big)N(x)=0 \ \ N(b)=1,$$
 $$
  1- (2/\rho) N^{\prime}(b)= \gamma \int_{0}^{b}e^{\gamma(z-b)} N(z)  d z, \  2\rho N^{\prime}(b)- 2N^{\prime\prime}(b)=\displaystyle{ \rho\gamma - \rho\gamma^{2}  \int_{0}^{b}e^{\gamma(z-b)} N(z)  d z   }$$ 
 
 Inserting $    N(x)=\alpha+\beta_- e^{s_- x}+\beta_+ e^{s_+ x}$ where $2s_{\pm}=\rho\pm\sqrt{\rho^2+4\gamma^2} $ results in a   linear system   for $\alpha, 
\beta_{\pm} $. 
After a  considerable  amount of algebra  the EP  simplifies to 

\beq N_b(x)=   \frac{A(x,b)}{A(b,b)}\label{ asym} \text{ where }$$
$$   A(x,b)=    s_+  e^{ (s_+-s_-)b}   \Big(1- \frac{s_-+\gamma  }{ \gamma} e^{  s_-x}  \Big) - s_-    \Big(1- \frac{s_++\gamma  }{ \gamma} e^{  s_+x}  \Big)  \enq  
Note how, despite being a Levy process,   the EP does not admit scale functions. 

 \end{exmp}

The following result gives the EP.   We skip the proof.  
 \begin{theorem}
 \label{Four4}       Suppose assumptions A1-A7 hold with $\tau_1\sim\mathcal E(\lambda)$ and    $h\in\mathcal H$  satisfies  \eqref{Four}.   Let $s_k, k=0,\dots n-1$ {\it be  the roots of Lundberg equation} 
   $      0=  (1-s/ \rho)Q(s)-R(s)$    and suppose they are all simple (Note that $s=0$ is always a root). Define  the $n\times n$ matrices  $ \mathbf{A}=(\mathbf{a}_{jk}) $  and $ \mathbf{M}=( m_{jk}),\  j, k=0\dots n$  with entries 
$$ m_{jk}:= \int_{-b}^0  \p_z^j h_- (y) e^{s_ky}dy  , \ \mathbf{a}_{0k}=1,\     k=0\dots n-1,$$ 
$$\mathbf{a}_{jk}=   s_k^{j-1} - s_k^{j} /\rho+(-1)^{j}m_{j-1,k}  + \sum_{l=0}^{j-2}(-1)^{j-l-1}I_{j-l-2}s_k^l, \  \ j=1,\dots n-1 $$ 
  Then the   EP is given by  \eqref{escape6} where
 $x'\eq x-b$ and  ${\bm \Theta}(x,b)$ is the $ (n+1)\times (n+1)$ matrix

\beq  {\bm \Theta} (x,b)= \left(\begin{matrix}0   &  e^{s_1 x'}&\dots & e^{s_{n-1}  x'}
 \\   1 &   1   &\dots &  1
 \\ 
   \bar H(0^+) &        \mathbf{a}_{10}    &        \dots& \mathbf{a}_{1,n-1}    \\    \vdots& &&  \vdots \\    (-1)^{n}h^{(n-1)}_+(0)  &   \mathbf{a}_{n-1,0}    &        \dots& \mathbf{a}_{n-1,n-1}  \end{matrix}\right)\label{Theta3} \enq 
 \end{theorem}


\subsection{The case of zero drift: $c=0$}   The drift-less  case   deserves particular interest for its relevance    to reliability theory.
Besides several  interesting  simplifications occur.   We reformulate Corollary 1 as

\begin{Corollary} Suppose  assumptions A1-A4 hold and that  $c=0$. Then the escape probability is independent of the history and of the arrival distribution $F$.  $N$  solves the integral equation  \eqref{26}
  
  \end{Corollary}

We now turn our attention to solving this under appropriate restrictions.  If either   $q_1=0$ or $p_1=0$ vanish   the solution becomes quite simple. 

 \begin{theorem} \label{gentrivialc=0} Suppose   that $c=0$ and let   $p= \Bbb P(J_1>0), q=1-p$.   \begin{enumerate}
\item $\vec p=(\star,\star,0,\star)$: Suppose 
support  $H \subset(-\infty,0)\cup[b-x,\infty)$     and   that $ -J_1\mathbf{1}_{J_1<0} $  
has density $h\eq h_-$  (see \eqref{H}).    The EP is 
   \beq   N_b(x) =     \frac{1}{2\pi i}     \int_{\beta-i\infty}^{\beta+i\infty}    \frac{p  e^{ xs} ds }{s\Big(1-   q\hat h_-(s)    \Big)}, 0\le x\le b \label{c=001}  \enq 
\item  $\vec p=(\star,0,\star,\star)$: Suppose 
support  $H \subset(-\infty,-b)\cup(0,\infty)$     and   that $  J_1\mathbf{1}_{J_1>0} $  
has density $h$. The EP is 
   \beq   N_b(x) =1-   \frac{1}{2\pi i}     \int_{\beta-i\infty}^{\beta+i\infty}    \frac{q  e^{(b-x)s} ds }{s\Big(1-  p \hat h_+(s)    \Big)} \label{c=00} \quad  \enq  
   \end{enumerate}
\end{theorem}  Note how this agrees with \eqref{reverse}.

{\it Proof}. The result follows   taking a LT on    \eqref{26}, which  reads, respectively  \beq    N  (x)=   p+ q  \int_{0}^{x}  N(x-
 z)   h_-(z)  dz \text{ and }$$$$   N  (x)=   \bar H(b-x)+ p  \int_{0}^{b-x}  N(x-
 z)   h_+(z)  dz \enq


We now consider the case when  severities have   rational   CF. The result follows from those of last section  letting $1/\rho=0$. 
  \begin{proposition} \label{ration}    Suppose $c = 0$ and   that  severities have   rational   CF given by \eqref{Four}. 
   \begin{enumerate}
\item  Let   $L \eq (Q -R)( s)$ and $\Bbb L\eq  L(\p_x)$. Then  $N(x)$  solves the  ODE  \eqref{solvable2} $\Bbb LN=0$  and the BC's \eqref{BC3} setting     $1/\rho\equiv 0$ 
  
 \end{enumerate}
\end{proposition}

\subsubsection{ Different examples}
\begin{exmp} \rm We consider the case      support   $J_1 =(-\infty,-b]\cup \{y_1\} $, namely   positive jumps
have a {\it  fixed magnitude } $y_1>0$  and negative jumps exceed $b$.     Let  $ p=\Bbb P(J_1=y_1)\in(0,1)$.  Recalling that $   k_1\eq  \lfloor (b-x)/y_1\rfloor$ and letting $c=0$   \eqref{solfix}  simplifies to 
 
 
 \beq  N_b(x)= 1- (1-p) \sum_{k=0}^{    k_1 } p^k  = p^{k_1+1}, x<b      \enq     
 Alternatively, note that  $(X)$ escapes through  $b$ iff  the first $k_1+1$ jumps  are  positive.
\end{exmp}
 
\begin{exmp}\rm   We study EP for the family of      CFs 
 \beq   \tilde h(\omega)=pe^{i\epsilon_+\omega  b} \frac{\gamma_+}{\gamma_+ +i\omega } +qe^{i\epsilon_-\omega  b} \frac{\gamma_- }{\gamma_- -i\omega } \label{five-parameter}  \enq 
depending on five  parameters $0\le p\le 1, q=1-p, \gamma_{\pm}>0$ and  
$\epsilon_{\pm}=0,\pm 1$ . We denote $$N(x):=  N_{(p,\gamma_+,\gamma_-,\epsilon_+,\epsilon_-)}(x) \text{  and }
   h(x)= p\gamma_+e^{-\gamma_+(x- \epsilon_{+} b)}\theta (x- \epsilon_{+} b)+ q\gamma_-e^{\gamma_-(x- \epsilon_{-} b)}\theta (-x+ \epsilon_{-} b))  $$ 
When  $   -\epsilon_{-}=  \epsilon_{+}=1 $ the problem is trivial:     $N(x)=p$. The case    $\epsilon_{-}=0,  \epsilon_{+}=1$  is covered by   theorem \eqref{gentrivialc=0}: the  EP   follows from  
\eqref{c=001}
\beq   
N_{(p,1,0,\gamma_+,\gamma_-)}(x)=1-q e^{-p\gamma_- x} \text{   }  \enq  
If 
  $   \epsilon_{-}=-1, \epsilon_{+}=0$  the EP     follows from  
\eqref{c=00}.   We find (which agrees with \eqref{reverse}) 
  \beq   N_{(p,0,-1,\gamma_+,\gamma_-)}(x)=1-N_{(q,1,0,\gamma_-,\gamma_+)}(b-x)= p e^{-q\gamma_+ (b-x)} \enq

     When  $\epsilon_{\pm}=0$ theorem \eqref{gentrivialc=0}  does not apply; nevertheless, since  it  corresponds to $h\in\mathcal H$, viz. \eqref{three-parameter}   the ideas of this section do.   It follows from  \eqref{three-parameterf}  and  \eqref{three-parameterR} that $$\bar H(0^+)=p, I_0=   p\gamma_+ -q\gamma_-,  v\eq  q\gamma_+- p\gamma_-$$  Here  $  Q-R(s)= -s^2+vs$. Hence the EP  can be found solving Eq. \eqref{solvable2}:   
  \beq \Bbb LN\eq  \Big(    \p_{xx}-v \p_x\Big)N(x)=0,   \text{ or } N(x)=\alpha + \beta e^{ v x}\label{solvable3}  \enq 
 with   appropriate  BC's \eqref{BC3}  with $n=2$ and     $1/\rho\equiv 0$. After tedious algebra one finds that the EP under jump density \eqref{three-parameter} $ h(x)= p\gamma_+e^{-\gamma_+ x }\theta (x )+ q\gamma_-e^{\gamma_-x }\theta (-x)$  is 


 \beq N_{(p,0,0,\gamma_+,\gamma_-)} (x)=p \frac{\gamma_--  q( \gamma_++\gamma_- ) e^{v x}    }{ p\gamma_-- q\gamma_+  e^{v b}  }  \label{c=0,asym} \enq

  \end{exmp}

  \begin{exmp} \rm   We consider different special cases 
\begin{enumerate} 

\item $p=\frac{\gamma_-}{\gamma_-+\gamma_+} $ or $\tilde h(\omega)= \frac{ \gamma_- \gamma_+ }{ (\gamma_+ -i\omega)( \gamma_- +i\omega)} $.  This gives  {\it variance gamma distribution (VGD) }  \eqref{variance}  with parameters      $n=1, \ \s^2=2\gamma_-\gamma_+   
\text{ and  }  \vartheta= \frac{\gamma_--\gamma_-}{\gamma_-+\gamma_+} $.   

 \item $p=\frac{\gamma_+}{\gamma_-+\gamma_+} $. This is a limit case for NPC: here   $\Bbb E J_1=v =0$.          
 Hence     the EP   \eqref{c=0,asym} is ill-defined   and  must be   obtained from scratch; we find 
 \beq N_b(x)=  
  \gamma_+  \frac{    \gamma_- x +1   }{      \gamma_++ \gamma_-+  b    
\gamma_+\gamma_-} \label{linear2} \enq 
  \item           Letting  $\gamma_+=\gamma_- $ in eq.  \eqref{linear2} one recovers   $ N(x)= \Big(1 +  \gamma x \Big)/\Big(2  +b  \gamma\Big)$, the EP under   $J_1 \sim$    Laplace$(0,\gamma)  $. Note how \eqref{symmetric} is      satisfied.   
  
\end{enumerate}
 
  \end{exmp}
\begin{rem}\rm  Notice   that the EP has the neat factorization $N_b(x)\eq  A(x)/B(b)$. This is interesting, since  $(X)$ is a two-sided not necessarily  L\'evy-Markov process.  
\rm

  Letting    $\gamma_{\pm} \to\infty $ with   $p=1/2$  and $2v\eq   \gamma_+- \gamma_- $ {\it constant  \eqref{c=0,asym} goes into   the scale function for Brownian motion} with drift $v$.     We have 
\beq  N_b(x)=\frac{\gamma_--  (\gamma_++\gamma_-)/2 e^{v x}    }{ \gamma_-- \gamma_+  e^{v b}  }  \underset{\gamma_{\pm} \to\infty}\to  \Big(1-   e^{ v x}     \Big)/\Big( 1-   e^{ v  b}   \Big) \enq 
   Actually   the infinitesimal generator of       BM is $\Bbb L\eq (1/2)\p_{xx}-v \p_x$ and the    EP   {\it satisfies also} \eqref{solvable3} $\Bbb L N=0$  with  {\it different   BCs } $N(b)=1-N(0)=1$.    Such remarkable coincidence can be traced to the   fact that if $ \tilde h( \omega \Delta)-1=O( |\omega \Delta|^\alpha), \Delta \to 0$ then marginal probabilities for \eqref{process} converge in a weak sense into a fractional diffusion, in particular to BM.    \end{rem}  
 



 \appendix 

 \section{}   \begin{theorem}   Suppose that  assumptions 5-7  hold  and let  $c>0, b<\infty$. 
 Then   \eqref{Fred} (or \eqref{fav2})  has a unique continuous solution  $x\mapsto N_b(x)$ which 
   satisfies the bound
 \begin{equation}0<        \sup_{x\in (0,b)} N(x) \le  \frac{L}{1-L} \label{bound}. \end{equation}  
where we recall that $ L:= \Bbb P\Big(\tau_1\le b/c, J_1\in (-b,b)\Big)  $, see   assumption A6.
Further,  \eqref{Fred} without forcing term has only the trivial solution.  
 \end{theorem}

{\it Proof}. We  use Banach fixed point theorem  with  the metric on  $  C_0(0,b) $  induced by the sup-norm  
on  $(0,b)$ denoted as   $\| . . \|_\infty$.   Given $ N\in C_0$  let us   introduce  the  integral operator  
 \beq   N   \mapsto \mathcal   O N(x):=  N^H(t_b)+  \tilde  {\mathcal    O} N(x) \text{ where } $$$$\ N^H(y)=  \bar F (y) + p_2  F (y) +  p_1 \int_0^{y}  dF(l)   \bar     H_{1+} ( c(y- l) ),   \enq 
  \beq  \tilde  {\mathcal    O} N(x)\eq 
  \int_0^{t_b}  dF(l)\Big(   
   \int_{0}^{b-x-cl}      p_1 d       H_{1+}(y) +    \int_{-x-cl}^0   q_1dH_{1-}(y)\Big)N(x+cl+y)   \enq      Let $N\in B\subset  C_0(0,b) $, the unit ball in $ C_0(0,b) $. Note first  that 
   $$ \|  \tilde  {\mathcal    O}  N\|_\infty\eq \underset   {x\in (0,b)} \sup \big|  \tilde  {\mathcal    O}  N(x)\big|  \le \|   N\|_\infty   \underset    {x\in (0,b)}  \sup
    \int_0^{t_b}  dF(l)\Big(   
   \int_{0}^{b-x-cl}      p_1 d       H_{1+}(y) +    \int_{-x-cl}^0   q_1dH_{1-}(y)\Big) =$$  
$$ \|   N\|_\infty   \underset    {x\in (0,b)}  \sup
    F(t_b) \Big(     p_1 d       H_{1+}(b-x) +       q_1 H_{1-}(-x)\Big) \le   L<1   $$
 Setting for convenience  $q_1=0, y:= b-x-cl$ and  recalling  assumption A5      we also  have 
$$    |    \tilde  {\mathcal    O} N(x+\e)-  \tilde  {\mathcal    O}  N(x) | \le  p_1\|N\|_\infty   \int_0^{t_b}  dF(l)   \big|  H_{1+}( y)-H_{1+}(y-\e)\big|\le     p_1\|N\|_\infty       F(t_b) \underset   {z\in (0,b)} \sup H_{1+}( z)-H_{1+}(z-\e), $$
$$\underset   { \e\to 0 } \lim \underset   {x\in (0,b)} \sup  |    \tilde  {\mathcal    O}N(x+\e)-  \tilde  {\mathcal    O}  N(x) | = p_1 F(t_b)  \underset   { \e\to 0 } \lim \underset   {z\in (0,b)} \sup H_{1+}( z)-H_{1+}(z-\e) =0$$  Alternatively, use Scheffe's
 theorem. Therefore   $  \tilde  {\mathcal    O}  B\subset  B  $:$  \tilde  {\mathcal    O} $ is a bounded endomorphism  of $B$.

More generally for $N_1, N_2\in C_0$ linearity implies  
  $$ \|  \mathcal    O   N_1-\mathcal   O  N_2 \|_\infty = 
\|   \tilde  {\mathcal    O}N_1 - \tilde  {\mathcal    O}N_2 \|_\infty= \|   \tilde  {\mathcal    O}(N_1 - N_2) \|_\infty  \le L \|   N_1-       N_2 \|_\infty$$  
Hence,  {\it  if   Assumptions 5,  6  hold   $\mathcal    O$  is    a contraction operator on $ C_0(0,b) $}:  Lipschitz continuous with Lipschitz constant $  L<1$.   
      Thus  it has just  one fixed point which satisfies $N(x)=\mathcal    O  N(x):= N^H(t_b)+  \tilde  {\mathcal    O} N(x) $, namely \eqref{fav2} and 
 $$\|  N\|_\infty \le \|  N^H(t_b)+  \tilde  {\mathcal    O} N(x) \|_\infty \le  \|  N^H \|_\infty+ L \| N\|_\infty$$
By  Gronwall's Lemma,  \eqref{fav2} without forcing term can only have the trivial solution  $N=0$. 
 
 \begin{rem} \rm 

When either $b\to\infty$ or  $c =q_2=p_2=0$   the Lipschitz constant $L\to  1$, the operator $\mathcal O$ is {\it not contractive, only  non-expansive}   and    \eqref{bound} blows up.   Actually the corresponding IE for survival probabilities {\it has always    a multi-parameter family of solutions}.

If, by contrast, Assumption A6 holds but A5 is dropped   the previous reasoning shows  mutatis-mutandis  that  there exist unique   solution but needs not being continuous. Finally, note that    Kolmogorov-Riesz   theorem  proves that $ \tilde  {\mathcal    O}$ is a {\it compact operator } on $C_0(0,b)$.

 \end{rem} 

 \section{}   {\it Proof of lemma \eqref{Lemma1}}. 
 Suppose         $f$  satisfies  \eqref{Lapl} and \eqref{Lapl2} where $Q,R\ne 0$. 
    Choose   the initial values $( f^{(k)}_0)$   to  satisfy the  system \eqref{linear}- which is  possible  since the associated system is  triangular     with  determinant $=a_n^n\ne0$.   Note that in this case   $\mathcal L\Big(Q(\p_t)f\Big)\eq  $
$$ \int_0^\infty e^{-st}    \Big(\underset{j=0}{\overset{ n}\sum}a_j  \f {\p^j \ }{\p t^j}\Big) f  (t)dt =     Q(s) \hat f(s) -    \underset{j=0}{\overset{ n-1}\sum} \Big(\sum_{k=0}^{n-j-1} a_{j+k+1}    f^{(k)}_0\Big)s^j =$$
$$    R(s)   -    \underset{j=0}{\overset{ n-1}\sum} \Big(\sum_{k=0}^{n-j-1} a_{j+k+1}    f^{(k)}_0\Big)s^j=   \underset{j=0}{\overset{ n-1}\sum} \Big(b_j-\sum_{k=0}^{n-j-1} a_{j+k+1}    f^{(k)}_0\Big)s^j=0$$ 
 By uniqueness of Laplace transform is $Q(\p_t)f=0$, i.e.  $f$ solves \eqref{Geneq} and \eqref{linear}.

  Reciprocally, suppose   that   $f$ solves \eqref{Geneq}  for  a certain minimal  $Q$ with IC  $( f^{(k)}_0)$; a simple calculation shows that   $\hat f(s)$ is  given by  \eqref{Lapl}   where  $P(s) = \underset{j=0}{\overset{ n-1}\sum}b_j  s^j $ and $b_j\eq \sum_{k=0}^{n-j-1} a_{j+k+1}    f^{(k)}_0, j=0,1,\dots n-1$. In particular $b_0=a_0$ and $\hat f(0)=1$.   $\ \square$

{\it Proof of Lemma 2}.   Since  $f$ is of class $C^n$ then   $   f ^{(j)}$  is bounded   and $ \underset{t\to 0}\lim f^{(j)}(t) =  f ^{(j)}_0=0$ exists, $j\le n$.  If  $  f ^{(j)}_0=0$ for some $1\le j\le n-1$ the  initial value theorem for LT    yields    
$$    0= \underset{t\to 0}\lim     f ^{(j)}   (t)=     \underset{s\to \infty }\lim       s^{j+1}\hat f(s)=    \underset{s\to \infty }\lim       s^{k+1}\hat f(s)=\underset{t\to 0}\lim     f ^{(k)}   (t), \forall  k\le j$$
 Besides $b_k\eq  \sum_{l=0}^{n-k-1} a_{l+k+1}    f^{l}_0   $ vanishes 
provided $n-k-1\le j\Leftrightarrow k\ge n-j-1$.


\begin{thebibliography}{99}


 \bibitem{Feller1} W. Feller 1954, Diffusion processes in one dimension.\textit{ Trans. Amer. Math. Soc.},  \textbf{77}, 1-3 

 



\bibitem{Bhat} R.N. Bhattacharya and  E.  C.Waymire 1981,  {\it Stoch. processes with Appl.}, Wiley, New York


\bibitem{kt81} S. Karlin and H. Taylor 1981. {\it A  first course in stochastic processes}, Acad. press, New York


       


 
  



 



\bibitem{Sornette}A. Helmstetter and D. Sornette 2003, Diffusion of epicenters of earthquake aftershocks, Omori law and generalized continuous-time random walk models. \textit{Phys. Rev. E}, \textbf{66}, 061104 

\bibitem{Merton} R. C.  Merton, Option pricing when   stock returns are discontinuous.  \textit{J.  Fin. Econ.}   3, 125-144

 

\bibitem{Cramer} H. Cramer 1994, On the mathematical theory of risk. \textit{Collected Works
vol. 1 (Springer)}, 601-678


  \bibitem{Sparre} Sparre E. Andersen 1957, 
 \textit{Transactions XVth Intern. Cong.
Actuar. New York, II}, 219-229,   





 

 \bibitem{Mikosch}T. Mikosch 2006,{  \it Non-Life Insurance Mathematics.}   Springer Verlag, New York 


 

 \bibitem{Rolski}T. Rolski, H.Schmidli, V.Schmidt, J.Teugels 2006, {\it Stochastic Processes for Insurance and Finance,} Wiley


\bibitem{dh0}  D. C. M. Dickson and C. Hipp 1998,   a class of renewal risk process. \textit{ North Amer.  Act. J.}, \textbf {7}, 1-12  
\bibitem{dh01} D. C.   Dickson, C. Hipp 1998,  Ruin probabilities for Erlang(2) 
process. \textit{ Insur. Math. Econ.}\textbf{22}, 251 
\bibitem{Garcia} J. M. Garcia 2005,  Explicit  
 survival probabilities in the classical risk model. \textit{ ASTIN Bull.}, \textbf {35}, 113-130 



 
\bibitem{Kluppelberg2} C. Kluppelberg 1989, Estimation of ruin probabilities by means of hazard rates. \textit{ Insur. Math.Econ.}, \textbf{8},  279-285 


 

 


  \bibitem{dh012} D. C.   Dickson 2001, C. Hipp, Time to ruin for Erlang(2)  
 process. \textit{Insur. Math. Econ.}, \textbf {29}, 333 




 \bibitem{Berger2}H.U. Gerber  and E.S.W. Shiu 2005, The time value of ruin in a Sparre Andersen model. \textit{N. Amer. Act. J.}, \textbf {9(2)}, 49-69 
 





  

\bibitem{lg04}  S. Li and J. Garrido 2004, On ruin for the Erlang(n) risk process. \textit{Insur. Math. Econ.}  \textbf {34}, 391-408 
 


 





\bibitem{dh08} D. C. M. Dickson, B.D. Hughes, Z. Lianzeng 2005,  The density of the time to ruin for   Sparre Andersen process with Erlang arrivals and exponential claims. \textit{Scand. Actuar. J.}, \textbf {5}, 358-376 

  




 



 



\bibitem{Kluppelberg3} C. Kluppelberg,  T Mikosch 1995,  Delay in claim settlement and ruin probability approximations. \textit{ Scand. Actuar.J.} \textbf{1995(2)}154-168 
 
 
\bibitem{Mikosch2}T. Mikosch, G. Samorodnitsky 2000, Ruin probability with claims modeled by a stationary ergodic stable process. {\it Ann.   Prob.} \textbf{28(4)} 1814
-1851 

\bibitem{Huzak}M. Huzak,  M. Perman, S. Hrvoje and Z. Vondracek 2004, Ruin probabilities
and decompositions for general perturbed risk processes. \textit{Ann.  Appl. Prob.}, \textbf{14},  1378-1397 
 \bibitem{Kluppelberg} C. Kluppelberg, A.E. Kiprianou, R.A. Maller 2004, Ruin probabilities and overshoots for general L\'evy insurance risk processes. \textit{Ann.Appl.Prob.}, \textbf{14(4)}, 1766-1801 



\bibitem{Wang} W Rongming, LHaifeng 2002,  Ruin Probability Under 
Risk  
Processes. \textit{ASTIN Bull}, \textbf{32},   81-90 

 
\bibitem{Dang} L Dang, N Zhub, H Zhangb 2009, Survival probability for 
 risk model. \textit{Insur.Math.Econ.}, \textbf{44(3)}, 491-496 



\bibitem{Lefevre}C. Lef\'evre,  S. Loisel 2007, 
Ruin probabilities for classical risk models. \textit{Scand. Actuar.J.}, 41-60 



 \bibitem{MV3} J. Villarroel, M. Montero 2011, Poisson driven stochastic    Schrodinger equation. \textit{Stud. Appl. Math.}\textbf{127}   372


 \bibitem{Gole}  L. Lavergnat,  P. Gole 1998,  Stochastic raindrop time distribution model. \textit{J. Appl. Meteor.}\textbf {37},  805


  
\bibitem{Perona2} P. Perona, E Daly, B Crouzy, A. Porporato 2012,  Stochastic dynamics of snow avalanche occurrence by superposition of Poisson processes. \textit{ Proc. R. Soc.A}, \textbf{468}, 4193-4208 

\bibitem{Othmer}H.  Othmer, S.  Dunbar, W. Alt,   Models of dispersal in biological systems 1988\textit{J. Math. Biol.}\textbf{26(3)}




 \bibitem{Oks} B.K. Oksendal 2005.{ \it Applied Stochastic Control of Jump Diffusions.} Springer Verlag, New York 

 

 \bibitem{Duan} J. Duan 2015, {\it An Introduction to Stochastic Dynamics 2,}  Cambridge University Press, New York 

 

 \bibitem{Liao} M. Liao 1989, The Dirichlet problem of a discontinuous Markov process. \textit{Acta Math. Sin.}, \textbf{5},   915 


\bibitem{Bertoin}  J. Bertoin 1996,  On the first exit time of a completely asymmetric stable process from a finite interval  \textit{  Bull. Lond. Math. Soc. }28, 514–520  


\bibitem{MV} M. Montero and J. Villarroel 2010, Mean exit times in non-Markovian drifting random-walk  processes. \textit{ Phys.  Rev. E:  Statistical Phys.}, \textbf {82}, 021102 

 \bibitem{Surya}  B. A. Surya 2008,  
Scale Functions of Spectrally Negative L\'evy   Processes
 \textit{  J. Appl. Prob. } {\bf 45},  135

\bibitem{Chiu}  S N Chiu and C Yin 2005,    Passage  times  for  a  spectrally  negative L\'evy  process  with  applications  to  risk theory
   \textit{  Bernoulli } \textbf {11(3)},  511-522 
\bibitem{Avram}  F. Avram, A. E. Kyprianou and M. R. Pistorius 2004,   Exit Problems for Spectrally Negative   L\'evy  Processes and Applications to   Russian Options,   \textit{ Ann. Appl. Prob. } {\bf 14(1)},  215-238

%
 






 





  


 




 




  







  





\bibitem{Feller3} W. Feller 1991,   Introduction to Probability Theory and   Applications, Vol. 2. \textit{Wiley} %

 

\bibitem{Sumita} U Sumita, Y Masuda 1987,  Classes of Probability Density Functions Having Laplace Transforms with Negative Zeros and Poles. \textit{ Adv. Appl. Prob.}, \textbf {19(3)},  632-651  


\vspace{0.3cm}
 
 
  




   



\end{thebibliography}
\end{document}